\documentclass[11pt,leqno]{amsart}

\usepackage{amsmath,amsthm}%,amsfonts}
\usepackage {latexsym}
\usepackage{amssymb}
\newcommand{\bse}{\begin{equation}}

\newcommand{\bear}{\begin{eqnarray}}
\newcommand{\eear}{\end{eqnarray}}

\newcommand{\supp}{\mbox{\rm supp}}

\newcommand{\eps}{{\varepsilon}}
\newcommand{\R}{{\mathbb R}}

\newcommand{\Z}{{\mathbb Z}}

\newcommand{\calS}{{\mathcal S}}

\newcommand{\shell}{{\mathcal S}}

\newcommand{\dist}{\mbox{\rm dist}}

\newcommand{\less}{\lesssim}
\newcommand{\les}{\lesssim}

\newcommand{\la}{\lambda}

\renewcommand{\phi}{\varphi}

\def\nn{\nonumber}

\newtheorem{theorem}{Theorem}
\newtheorem{lemma}[theorem]{Lemma}

\newtheorem{cor}[theorem]{Corollary}
\newtheorem{prop}[theorem]{Proposition}

\theoremstyle{remark}
\newtheorem{remark}[theorem]{Remark}

\def\calJ{{\mathcal J}}

\def\la{\langle}
\def\ra{\rangle}

\def\norm[#1][#2]{\|#1\|_{#2}}

\def\x{{\bf x}}

\def\norm[#1][#2]{\Vert #1 \Vert_{#2}}

\renewcommand{\div}{{\rm div}}
\def\wenn{\text{\ \ \rm if\ \ }}
\def\PV{{\rm P.V.}}
\def\chitil{\tilde{\chi}}
\def\FT{{\mathcal F}}
\renewcommand{\hat}{\widehat}
\renewcommand{\bar}{\overline}
\def\pj{{(j)}}
\def\pone{{(1)}}
\def\ptwo{{(2)}}
\renewcommand\tilde{\widetilde}
\begin{document}

\title[Strichartz estimates for large magnetic potentials]
{Strichartz and smoothing estimates for Schr\"odinger operators with
large magnetic potentials in~$\R^3$}
\author{M. Burak Erdo\smash{\u{g}}an, Michael Goldberg, Wilhelm Schlag}

\address{Department of Mathematics \\
University of Illinois \\
Urbana, IL 61801, U.S.A.}
 \email{berdogan@math.uiuc.edu}

\address{Department of Mathematics\\ Johns Hopkins University\\
Baltimore, MD 21218, U.S.A.}

\email{mikeg@math.jhu.edu}

\address{ Department of Mathematics, University of Chicago,
5734 South University Avenue, Chicago, IL 60637, U.S.A.}
\email{schlag@math.uchicago.edu}

\maketitle

\begin{abstract}
We show that the time evolution of the operator
\[
H= -\Delta+ i(A\cdot\nabla + \nabla\cdot A) + V
\]
in $\R^3$ satisfies global Strichartz and smoothing estimates under
suitable smoothness and decay assumptions on $A$ and $V$ but without
any smallness assumptions. We require that zero energy is neither an
eigenvalue nor a resonance.
\end{abstract}

\section{Introduction}

Magnetic Schr\"odinger operators on $L^2(\R^d)$ are of the form
\begin{equation}
\label{eq:schr1} H = -\Delta+ i(A\cdot\nabla + \nabla\cdot A) + V = -\Delta + L
\end{equation}
There has been much activity surrounding dispersive estimates for
the case $A=0$ under suitable decay (and also regularity when
$d\ge4$) assumptions on~$V$. In fact, in that case the harder
$L^1(\R^d)\to L^\infty(\R^d)$ estimate is now known in all
dimensions $d\ge1$ under the condition that zero energy is neither
an eigenvalue nor a resonance (and there are now also results in the
case when the latter assumption does not hold). The seminal paper
for this class of estimates is~\cite{JSS} and we refer the reader
to~\cite{Sch} for a survey of more recent work.

On the other hand, much less is known when $A\not\equiv0$. In
\cite{Stef} and~\cite{GST} Strichartz and smoothing estimates were
obtained for small $A$ and $V$. In this paper we prove the following
theorem:

\begin{theorem}
\label{thm:main} Let $A$ and $V$ be real-valued such that for all
$x,\xi\in\R^3$
  \begin{align}
  \la x\ra|A(x)|+|DA(x)|+|V(x)| &\less \la x\ra^{-8-\eps}  \label{eq:ass1}\\
\sum_{|\alpha|\le 2} |D^\alpha
  \hat{A}(\xi)| &\less \la \xi\ra^{-3-\eps}  \label{eq:ass2} \\
  |\nabla V(x)| &\less \la x\ra^{-1-\eps}  \label{eq:ass3}
\end{align}
for some $\eps>0$. Furthermore, assume that zero energy is neither
an eigenvalue nor a resonance of $H$. Then, with $P_c$ being the
projection onto the continuous spectrum,
\begin{equation}\label{eq:strich} \|e^{itH} P_c f\|_{L_t^q(L_x^p)}
\less \|f\|_{L^2(\R^3)} \end{equation} provided
$\frac{2}{q}+\frac{3}{p}=\frac{3}{2}$ and $2\le p<6$. Moreover, the
Kato smoothing estimate
\begin{equation}\label{eq:smooth}
\int_0^\infty \big\|\la x\ra^{-\sigma}\la \nabla\ra^{\frac12}
e^{itH}P_c f\big\|_2^2\, dt \le C\|f\|_2^2
\end{equation}
holds with $\sigma>4$.
\end{theorem}

The definition of zero energy being neither an eigenvalue nor a
resonance is the usual one:  there does not exist $f\in
\cap_{\tau>\frac12} L^{2,-\tau}(\R^3)$, $f\not\equiv0$ such that
$Hf=0$.

In a sequel to this paper the authors will weaken the conditions on
$A$ and $V$ --- in fact, for the sake of simplicity we have chosen
to impose somewhat stronger conditions on $A$ and $V$ than the
methods of this paper actually require. Let us merely comment
that~\eqref{eq:ass3} can be dispensed with -- we only include it in
order to keep this paper self-contained: it is used to prove the
absence of imbedded eigenvalues as in the appendix. However, in the
recent work~\cite{KochT} a much stronger result is presented that
does not require this condition.

 The approach in this work is perturbative around the free case despite the fact that we
make no smallness assumption; instead, we use Fredholm theory as
usual. The actual perturbation argument is the one from~\cite{RS}
where it was used in the case of $A=0$. The main novel ingredient in
this paper is a limiting absorption estimate for large energies.
More precisely,  recall that in~\cite{Agmon} and~\cite{IS} it is
proved that for $H$ as in~\eqref{eq:schr1} under suitable decay
conditions on $A$ and $V$ and with $\tau>\frac12$,
\begin{equation}\label{eq:limap0} \sup_{\lambda\in[\delta,\delta^{-1}]} \|
\la \nabla\ra \la x\ra^{-\tau} (H-(\lambda^2+i0))^{-1}\la
x\ra^{-\tau} \la \nabla\ra \|_{2\to2} \le C(\delta)<\infty
\end{equation}
provided there are no imbedded eigenvalues in the continuous
spectrum.  It is well--known that this {\em limiting absorption
principle} is of fundamental importance for proving dispersive
estimates, at least for the case of large potentials. However, for
this one needs to remove the restriction on $\lambda$. To
extend~\eqref{eq:limap0} to zero energies is similar to the case
$A=0$. This step requires the assumption on zero energy.

Note that~\eqref{eq:limap0} as stated cannot hold as
$\lambda\to\infty$ since it fails even for the free resolvent.
Indeed, with $\tau>\frac12$
\begin{equation}\label{eq:freeone}
\sup_{\lambda>1}\|\la \nabla\ra^{\frac12}\la x\ra^{-\tau}
(H_0-(\lambda^2+i0))^{-1}\la x\ra^{-\tau}\la
\nabla\ra^{\frac12}\|_{2\to2} \less 1
\end{equation}
and this is optimal in the sense that no more than one derivative in
total can be gained here.  We will adopt the shorthand notation
\[ R_0(z) := (H_0 - z)^{-1} \]
for the resolvent of the Laplacian.  The resolvent of a general operator
$H$ will be indicated by $R_H(z)$, or else $R_L(z)$ in the case where 
$H$ is specifically of the form $H_0 + L$.

In this paper we extend~\eqref{eq:freeone}
to~$H = H_0 + L$ for the class of first-order perturbations described in
Theorem~\ref{thm:main}.  A unified statement of the mapping properties
of the resolvent of $H$ over the entire spectrum $\lambda > 0 $ is as follows.
\begin{theorem} \label{thm:limap}
Suppose $H$ is a magnetic Schr\"odinger operator whose potentials satisfy the 
considitions \eqref{eq:ass1}--\eqref{eq:ass3}.  Then for $\tau > 4$
and $\alpha \in [0, 1]$, 
\begin{equation} \label{eq:limap1}
\sup_{\lambda > 1} \, \lambda^{1-2\alpha} 
\| \la\nabla\ra^{\alpha} \la x\ra^{-\tau}  (H - (\lambda^2 + i0))^{-1} 
   \la x \ra^{-\tau} \la\nabla\ra^{\alpha}\|_{2\to2} \les 1.
\end{equation}
If one further assumes that zero is not an eigenvalue or resonance of $H$,
then this bound can be extended to
\begin{equation} \label{eq:limap2}
\sup_{\lambda \ge 0} \ \la \lambda \ra^{1-2\alpha}
\| \la\nabla\ra^{\alpha} \la x\ra^{-\tau}  (H - (\lambda^2 + i0))^{-1}
   \la x \ra^{-\tau} \la\nabla\ra^{\alpha}\|_{2\to2} \les 1.
\end{equation}
As a consequence, the spectrum of $H$ is purely absolutely continuous over the
entire interval $[0,\infty)$.
\end{theorem}
\begin{remark} A result of type \eqref{eq:limap1}, in the case $\alpha=0$,
 is proved in \cite{Rob}
using the method of Mourre commutators.  In that work the potentials require
only very slight polynomial decay, however they are also
assumed to be infinitely differentiable, with the derivatives satisfying a 
symbol-like decay condition.
\end{remark}

Results of this type often rely upon the invertibility of the
operator $I + R_0(\lambda^2+i0)L$ in a suitable weighted space
$L^{2,-\sigma}$.   In the scalar ($A=0$) case, this becomes easy for
large $\lambda$ as the norm of $R_0(\lambda^2+i0)V$ decreases to zero.

One difficulty encountered here is that the norm of $R_0(\lambda^2 + i0)L$ 
does not decay as $\lambda \to \infty$,
since there is no decay to be found in the operator 
estimate~\eqref{eq:freeone}.
%Unlike the case $A=0$ the norm of
%$R_0(\lambda^2+i0)L$ does not become small as $\lambda\to\infty$
%since in~\eqref{eq:freeone} there is no decay of the left-hand side
%for increasing~$\lambda$. 
To circumvent this, we reduce ourselves to
the invertibility of $I-(-1)^m(R_0(\lambda^2+i0)L)^m$ and show that
$(R_0(\lambda^2+i0)L)^m$ is of small norm provided $m$ and $\lambda$
are large.

\section{The basic setup}

The following result is proved in \cite{RS}, see Theorem~4.1 in that
paper. It is based on Kato's notion of smoothing operators,
see~\cite{Kato}. We recall that for a self-adjoint operator $ H$, an
operator $\Gamma$ is called $H$-smooth in Kato's sense if for any
$f\in  {\mathcal D}(H_0)$
\begin{equation}
\label{eq:Csmooth} \|\Gamma e^{it H} f\|_{L^2_t L^2_x}\le
C_{\Gamma}(H) \|f\|_{L^2_x}
\end{equation}
or equivalently, for any $f\in L^2_x$
\begin{equation}
\label{eq:CsmoothR} \sup_{\eps >0} \|\Gamma R_{H}(\lambda\pm i\eps)
f\|_{L^2_\lambda L^2_x} \le C_{\Gamma}(H)\|f\|_{L^2_x}.
\end{equation}
We shall call $C_{\Gamma}(H)$ the smoothing bound of $\Gamma$
relative to $H$. Let $\Omega\subset \R$ and let $P_\Omega$ be a
spectral projection of $H$ associated with a set $\Omega$. We say
that  $\Gamma$ is $H$-smooth on $\Omega$ if $\Gamma P_{\Omega}$ is
$H$-smooth. We denote the corresponding smoothing bound by
$C_\Gamma(H,\Omega)$. It is not difficult to show (see
e.g.~\cite{RS4}) that, equivalently,  $\Gamma$ is $H$-smooth on
$\Omega$ if
\begin{equation}
\label{eq:smon} \sup_{\beta >0} \|\chi_\Omega(\lambda)\Gamma
R_{H}(\lambda\pm i\beta ) f\|_{L^2_\lambda L^2_x} \le C_{\Gamma}(H,
\Omega)\|f\|_{L^2_x}.
\end{equation}

The estimate \eqref{eq:strich} of Theorem~\ref{thm:main} is obtained
by means of the following result. The remainder of the paper is
devoted to verifying the conditions needed in
Proposition~\ref{prop:Katoth}. Furthermore, this verification will
establish the smoothing estimate~\eqref{eq:smooth}.

\begin{prop}\label{prop:Katoth}
Let $H_0=-\Delta$ and $H= H_0 + L$ with $L=\sum_{j=1}^J Y_j^* Z_j$.
We assume that each $Y_j$ is $H_0$ smooth with a smoothing bound
$C_B(H_0)$ and that for some $\Omega\subset \R$ the operators $Z_j$
are $H$-smooth on $\Omega$ with the smoothing bound $C_A(H,\Omega)$.
Assume also that the unitary semigroup $e^{it H_0}$ satisfies the
estimate
\begin{equation}
\label{eq:strH0} \|e^{it H_0} \psi_0\|_{L^q_t L^r_x} \le C_{H_0}
\|\psi_0\|_{L^2_x}
\end{equation}
for some $q\in (2,\infty]$ and $r\in [1,\infty]$. Then the semigroup
$e^{itH}$ associated with $H = H_0 + L$, restricted to the spectral
set $\Omega$, also verifies the estimate \eqref{eq:strH0}, i.e.,
\begin{equation}
\label{eq:strHk} \|e^{it H}P_\Omega \psi_0\|_{L^q_t L^r_x} \le J
C_{H_0}C_B(H_0) C_A(H,\Omega) \|\psi_0\|_{L^2_x}
\end{equation}
\end{prop}

We refer the reader to \cite{RS} for the proof. Note that this
approach does not capture the Keel-Tao endpoint (which would
correspond to $q=2$) --- the reason being the Christ-Kiselev
lemma~\cite{CK} which is used in the proof of
Proposition~\ref{prop:Katoth}. To apply this proposition we write,
with a decreasing weight $w(x)=\la x\ra^{-\sigma}$, for some
sufficiently large $\sigma>0$,
\begin{equation}\label{eq:Ldecomp}
\begin{aligned}
 L &= 2iA\cdot\nabla + i\div A + V \\
 &= 2i Aw^{-1} \cdot \nabla \la \nabla\ra^{-\frac12} \la
 \nabla\ra^{\frac12} w + 2i A\cdot \nabla (w^{-1}) w + i\div A + V
 \\
 &= \sum_{j=1}^2 Y_j^* Z_j
\end{aligned}
\end{equation}
where
\begin{equation}
\label{eq:YZdef}
\begin{aligned}
 & Y_1^* := 2i Aw^{-1} \cdot \nabla \la \nabla\ra^{-\frac12}, \quad Z_1
  := \la
 \nabla\ra^{\frac12} w \\
 & Y_2^* := [2i A\cdot \nabla (w^{-1}) w + i\div A + V] w^{-1}, \quad
  Z_2 := w
\end{aligned}
\end{equation}
Throughout this paper, we shall treat $\sigma>0$ as a parameter. In
various places we shall specify how large it needs to be chosen.
Eventually, we shall require $\sigma>4$, which will lead to the
condition~\eqref{eq:ass1}. It is standard that $Y_1$ and $Y_2$ are
$H_0$-smooth provided
\begin{equation}\label{eq:int_cond} |A(x)|+|\div A(x)|+|V(x)|\less \la
x\ra^{-1-\sigma-\eps}
\end{equation}
We now start discussing the smoothing properties of $Z_1$ and $Z_2$
relative $H$. It will suffice to discuss~$Z_1$.

Let us first consider intermediate energies $\lambda^2$, i.e.,
$\lambda\in[\lambda_0^{-1},\lambda_0]=\calJ_0$ with $\lambda_0$
large. Then it was shown in \cite{IS}, see also~\cite{Agmon}, that
the resolvent of $H$ satisfies the following bound
\[
\sup_{\lambda\in\calJ_0}\|\la x\ra^{-\frac12-\eps}\la \nabla\ra
R_L(\lambda^2+i0) f\|_2 \le C(\lambda_0)\, \|\la x\ra^{\frac12+\eps}
\la \nabla\ra^{-1} f\|_2
\]
(in fact, a stronger bound was proved in \cite{IS}). More precisely,
this bound follows provided there are no eigenvalues of $H$ in the
interval $\calJ_0$. However, we prove the latter property in the
appendix (see also~\cite{KochT}). Therefore,
\begin{equation}\label{eq:pseudo}
\begin{aligned}
\sup_{\lambda\in\calJ_0}\|Z_1 R_L(\lambda^2+i0) Z_1^* \|_{2\to2}
&\le C(\lambda_0) \|\la\nabla\ra^{\frac12}w \la \nabla\ra^{-1} \la
x\ra^{\frac12+\eps}\|^2_{2\to2} \le C(\lambda_0)
\end{aligned}
\end{equation}
since $\|\la\nabla\ra^{\frac12}w \la \nabla\ra^{-1} \la
x\ra^{\frac12+\eps}\|_{2\to2}<\infty$ by pseudo-differential
calculus. Finally, by Kato's smoothing theory, see \cite{RS4}
Theorem~XIII.30,  we conclude that $Z_1$ is $H$-smooth on
$\Omega=\calJ_0$.

Note that this argument does not carry over to $\lambda\to\infty$
(in other words, for magnetic potentials, unlike the case of $V$
alone, large energies are not easy). This is due to the fact that
the limiting absorption principles in \cite{IS} and \cite{Agmon} do
not yield a gain of one derivative uniformly in $\lambda$. We devote
Section~\ref{sec:large} to this issue.

Next, we turn to small energies.

\section{Small energies} \label{sec:small}

As usual, this is reduced to zero energy. For the latter, we need to
impose an invertibility condition which amounts to boundedness of
the resolvent $R_L(0)$ between suitable spaces. More precisely, by
the resolvent identity,
\[ R_L(\lambda^2+i0)= (1+R_0(\lambda^2+i0)L)^{-1}R_0(\lambda^2+i0)
\]
provided the inverse on the right-hand side exists. Therefore,
\begin{align*}
&\| Z_1 R_L(\lambda^2+i0) Z_1^* \|_{2\to2} \\
&=  \| Z_1
(1+R_0(\lambda^2+i0)L)^{-1} Z_1^{-1} Z_1 R_0(\lambda^2+i0) Z_1^*
 \|_{2\to2} \\
 &\le \| Z_1
(1+R_0(\lambda^2+i0)L)^{-1} Z_1^{-1} \|_{2\to2} \|Z_1
R_0(\lambda^2+i0) Z_1^*
 \|_{2\to2}
\end{align*}
By the smoothing properties of $Z_1$ relative to $H_0$,
\[
\sup_{\lambda}\|Z_1 R_0(\lambda^2+i0) Z_1^*
 \|_{2\to2} <\infty
\]
provided $\sigma>1$. For $\lambda>1$ this follows from
Agmon~\cite{Agmon} with $\sigma>\frac12$, whereas for small
$\lambda$ this can be reduced to a Hilbert-Schmidt norm provided
$\sigma>1$, see~\cite{JenKat}.

Thus, we need to verify that
\[
\begin{aligned}
&\sup_{|\lambda|<\lambda_0^{-1}}\| Z_1 (1+R_0(\lambda^2+i0)L)^{-1}
Z_1^{-1} \|_{2\to2} \\&=\sup_{|\lambda|<\lambda_0^{-1}}\| \la
 \nabla\ra^{\frac12} w
(1+R_0(\lambda^2+i0)L)^{-1} w^{-1}\la
 \nabla\ra^{-\frac12}  \|_{2\to2} <\infty
 \end{aligned}
\]
for some choice of large $\lambda_0$. First, we consider the case
$\lambda=0$. As usual, we let $G:=R_0(0)$.

\begin{lemma}
\label{lem:G} Assume that $L=2i\nabla\cdot A - i\div A + V$
satisfies $|A(x)|\less \la x\ra^{-\sigma-1-\eps}$, $|\div
A(x)|+|V(x)|\less \la x\ra^{-2\sigma}$ with $\sigma>1$. Then $Z_1 GL
Z_1^{-1}$ is a compact operator on $L^2$.
\end{lemma}
\begin{proof}
First, we consider only the $2i\nabla\cdot A$ part of $L$. We claim
that
\begin{equation}\label{eq:claim10} \|\la \nabla\ra G
\nabla\cdot A w^{-1} f\|_2 \less \|f\|_2
\end{equation}
To see this,
observe that by Plancherel
\[
\|D^\alpha G\nabla\cdot Aw^{-1} f\|_2 \less \| Aw^{-1} f\|_2 \less
\|f\|_2
\]
provided $|\alpha|=1$. On the other hand, we will show that
\begin{equation}\label{eq:Gnabla}
\|G\nabla\cdot Aw^{-1}f\|_2  \less \|Aw^{-1}f\|_{L^{2,1+\eps}} \less
\|f\|_2
\end{equation}
It suffices to prove that multiplication by $\frac{\xi}{|\xi|^2}$
maps $H^{1+\eps}$ to $L^2$. Let $\chi(\xi)$ be a smooth cut-off
around zero. Then $(1-\chi(\xi))\frac{\xi}{|\xi|^2}$ maps
$H^{1+\eps}$ to itself which is even stronger. Moreover, by
H\"older's inequality and Sobolev imbedding,
\[
\|\chi(\xi)|\xi|^{-1} g\|_2 \le \|\chi(\xi)|\xi|^{-1} \|_{L^{3-}}
\|g\|_{L^{6+}}\less \|g\|_{H^{1+\eps}}
\]
which implies \eqref{eq:Gnabla}. In conclusion, we have
proved~\eqref{eq:claim10}.

Thus,
\[\la
 \nabla\ra^{\frac12} w G \nabla\cdot A w^{-1} \la
 \nabla\ra^{-\frac12} = \la \nabla\ra^{\frac12} w \la \nabla\ra^{-1} \la \nabla\ra G \nabla\cdot A w^{-1}
 \la
 \nabla\ra^{-\frac12}
\] is compact in $L^2$, since
$\la\nabla\ra^\frac12 w \la\nabla\ra^{-1}$ is compact in $L^2$.

Second, we discuss the $\tilde V:= - i\div A + V$ part of $L$. It
will suffice to show that \begin{equation}\label{eq:lzwei} \| \la
\nabla\ra^{\frac12} w G\tilde V w^{-1}\la x\ra^{\eps} f\|_2\less
\|f\|_2
\end{equation}
since then
\[
\la \nabla\ra^{\frac12} w G\tilde V w^{-1} \la \nabla\ra^{-\frac12}
= \la \nabla\ra^{\frac12} w G\tilde V w^{-1} \la x\ra^{\eps} \la
x\ra^{-\eps}  \la \nabla\ra^{-\frac12}
\]
is compact. To prove \eqref{eq:lzwei}, we argue as before:
\[
\| \la \nabla\ra^{\frac12} w G\tilde V w^{-1}\la x\ra^{\eps}
f\|_2\less \| \nabla w G\tilde V w^{-1}\la x\ra^{\eps} f\|_2 + \| w
G\tilde V w^{-1}\la x\ra^{\eps} f\|_2
\]
The second summand on the right-hand side is controlled by the
Hilbert-Schmidt norm provided $\sigma>1$. The first summand is
similar to the proof of \eqref{eq:Gnabla}.
\end{proof}

The following remark will be used to analyze the condition at energy
zero.

\begin{remark}\label{rem:resonance}
Combining \eqref{eq:claim10} with the usual boundedness properties
of $G$ on weighted $L^2$ spaces (i.e., $G:L^{2,\beta_1}\to
L^{2,-\beta_2}$, provided $\beta_1+\beta_2>2$ and
$\beta_1,\beta_2>\frac12$, see \cite{JenKat} or~\cite{GS}) yields
\begin{equation}\label{eq:JenKat}
  \|GL h\|_{L^{2,-\tau+\eps/2}(\R^3)}\le
  \|h\|_{L^{2,-\tau}(\R^3)}
\end{equation}
for any $\tau>(1+\eps)/2$ provided $|\div A(x)|+|V(x)|\less \la
x\ra^{-2-\eps}$ and $|A(x)|\less \la x\ra^{-\tau-1-\eps}$.
\end{remark}

As an immediate consequence we arrive at the following.

\begin{cor}
  \label{cor:invert}
Assume that $\ker(I+Z_1 GL Z_1^{-1})=\{0\}$ as an operator on
$L^2(\R^3)$. Then $I+ Z_1 GL Z_1^{-1}$ is invertible on $L^2$.
Moreover,
\begin{equation}\label{eq:gleich}
\| Z_1 (I+R_0(\lambda^2+i0)L)^{-1} Z_1^{-1} \|_{2\to2} < \infty
\end{equation}
uniformly for small $\lambda$. An analogous statement holds with
$Z_2$ instead of $Z_1$.
\end{cor}
\begin{proof}
The first statement is Fredholm's alternative. Note that \[(I+ Z_1
GL Z_1^{-1} )^{-1} = Z_1  (I+GL)^{-1} Z_1^{-1}
\]
where $GL$ on the right-hand side is an operator on
$Z_1^{-1}(L^2(\R^3))$. By the same token, \eqref{eq:gleich} is the
same as
\[
\| (I+Z_1R_0(\lambda^2+i0)L Z_1^{-1})^{-1}  \|_{2\to2} < \infty
\]
uniformly for small $\lambda$.  To prove this, we write
\[
I+Z_1R_0(\lambda^2+i0)L Z_1^{-1} = I+ Z_1 GL Z_1^{-1}  + Z_1
B_\lambda L Z_1^{-1}
\]
where $B_\lambda = R_0(\lambda^2+i0) - G$. By a Neumann series
argument, it suffices to prove that
\begin{equation}
\sup_{|\lambda|<\lambda_0^{-1}} \|Z_1 B_\lambda L Z_1^{-1}\|_{2\to2}
\to 0
  \label{eq:claim3}
\end{equation}
as $\lambda_0\to\infty$. We have the following bounds on the kernel
of $B_\lambda(x,y)$:
\begin{equation}
\begin{aligned}
|B_\lambda(x,y)| &\less \frac{|\lambda|^\gamma}{|x-y|^{1-\gamma}},
\qquad 0\le \gamma\le 1 \\
|\nabla_x B_\lambda(x,y)\nabla_y| &\less \frac{\lambda}{|x-y|^2} +
\frac{\lambda^2}{|x-y|} \\
|\nabla_x B_\lambda(x,y)|+|B_\lambda(x,y)\nabla_y| &\less
\frac{\lambda}{|x-y|}
\end{aligned}
  \label{eq:Bbd}
\end{equation}
To prove \eqref{eq:claim3}, we estimate
\begin{align*}
\|Z_1 B_\lambda L Z_1^{-1}\|_{2\to2} &\less \|\nabla w B_\lambda L
w^{-1} \|_{2\to2} + \| w B_\lambda L w^{-1} \|_{2\to2} \\
&\less \|w \nabla B_\lambda L w^{-1} \|_{2\to2} + \| w B_\lambda L
w^{-1} \|_{2\to2}
\end{align*}
As before, we write $L=2i\nabla\cdot A + \tilde V$. To conclude the
argument, one now uses \eqref{eq:Bbd} together with Schur's lemma
(for the $\frac{\lambda}{|x-y|^2}$ term) as well as the
Hilbert-Schmidt norm (for the others).
\end{proof}

We now relate the condition in Corollary~\ref{cor:invert} to the
notion of resonance and/or eigenvalue at zero.

\begin{lemma}
  \label{lem:reson}
  Suppose that zero is neither an eigenvalue nor a resonance of
  $H$. Then under the conditions of Lemma~\ref{lem:G} one has
  \[
\ker(I+Z_j GL Z_j^{-1})=\{0\} \text{\ \ on\ \ }L^2(\R^3)
  \]
  for $j=1,2$. In particular, \eqref{eq:gleich} holds for small~$\lambda$.
\end{lemma}
\begin{proof}
Suppose $f\in L^2(\R^3)$ satisfies
\[
f+ Z_1 GL Z_1^{-1} f =0
\]
Set $h:= Z^{-1} f$. Then $h=-GL h\in L^{2,-\sigma}(\R^3)$. Applying
Remark~\ref{rem:resonance} we see that $h\in
L^{2,-(\sigma-\frac{\eps}{2})}(\R^3)$. Repeating this process shows
that $h\in \cap_{\tau>\frac12} L^{2,-\tau}(\R^3)$. It follows, see
\cite{JenKat} and~\cite{GS} that $Hh=0$ in the distributional sense.
However, by our assumption on zero energy it follows that $h=0$ and
therefore $f=0$ as desired. The argument for $Z_2$ is analogous.
\end{proof}

\section{Large energies}\label{sec:large}

The goal of this section is to prove the bound
\begin{equation}
  \label{eq:largelam}
  \sup_{\lambda>\lambda_0}\| Z_1 R_L(\lambda^2+i0) Z_1^* \|_{2\to2} <\infty
\end{equation}
with some large $\lambda_0$ and similarly with $Z_2$. Here $Z_1,
Z_2$ are as in \eqref{eq:YZdef} with $w(x)=\la x\ra^{-\sigma}$. Note
that in combination with the previous sections this will finish the
proof of Theorem~\ref{thm:main}. In order to
establish~\eqref{eq:largelam} we introduce some notations: for any
$\lambda>1$ define
\[
\widehat{T_\lambda f} (\xi) = \la \xi/\lambda\ra^{-1}\, \hat{f}(\xi)
\]
as well as
\[
S_\lambda := T_\lambda^{-1} R_0(\lambda^2+i0)
\]
It is clear that  for any $\tau$ one has
\begin{equation}\label{eq:Tlam}
T_\lambda : L^{2,\tau}\to L^{2,\tau} \end{equation} with a bound
independent of~$\lambda$. Indeed, by the Fourier transform this is
equivalent to
\[
\la \xi/\lambda\ra^{-1}\: :\: H^{\tau}\to H^{\tau}
\]
as a multiplication operator with norm independent of $\lambda$.
The decay in large $|\xi|$ suggests that $T_\lambda$ also improves local
regularity.  More precisely,
\[ \| \la\nabla\ra^\alpha T_\lambda f\|_{L^{2\tau}} \les \la \lambda\ra^\alpha
   \| f \|_{L^{2,\tau}} \]
for any $\alpha$ in the range $[0,1]$.

The Fourier multiplier associated to $S_\lambda$ is less well behaved,
however we still have the following bound:

\begin{lemma}
\label{lem:Sagmon} With $S_\lambda$ as before
\[ \| \la \nabla \ra^\alpha S_\lambda f\|_{L^{2,-\tau}} \less 
 \lambda^{\alpha-1} \|f\|_{L^{2,\tau}} \] provided $\tau>\frac12$ and
$\alpha \in [0,1]$.
\end{lemma}
\begin{proof}
By algebra of operators,
\begin{equation}\label{eq:algebra}
\la \nabla/\lambda\ra^2 R_0(\lambda^2+i0) =
2R_0(\lambda^2+i0)-\lambda^{-2}I
\end{equation}
Therefore, if $\tau>\frac12$ and $\lambda>1$, then
\begin{align*}
\|\la \nabla/\lambda\ra^2 R_0(\lambda^2+i0)f\|_{L^{2,-\tau}} &\le
2\|R_0(\lambda^2+i0)f\|_{L^{2,-\tau}}+\lambda^{-2}
\|f\|_{L^{2,-\tau}} \\
&\less \lambda^{-1} \|f\|_{L^{2,\tau}}
\end{align*}
by Agmon's limiting absorption principle~\cite{Agmon}. Finally, we
bound
\[
\| \la\nabla\ra^\alpha S_\lambda f\|_{L^{2,-\tau}}\le
\|\la\nabla\ra^\alpha T_\lambda\|_{L^{2,-\tau}\to L^{2,-\tau} } \|\la
\nabla/\lambda\ra^2 R_0(\lambda^2+i0)f\|_{L^{2,-\tau}}
\]
which finishes the proof.  
\end{proof}
\begin{remark}
The resolvent estimate that we used above, 
\[ \|R_0(\lambda^2+i0)f\|_{L^{2,-\tau}} \les \lambda^{-1} \|f\|_{L^{2,\tau}} \]
follows directly from the calculations in~\cite{Agmon}, but only appears as
a separately stated theorem in later works such as~\cite{JenKat}.
\end{remark}

Next, we combine $T_\lambda$ and $S_\lambda$ with $Z_1$ (in what
follows, we will treat $Z_1$, the case of $Z_2$ being easier):

\begin{lemma}
\label{lem:STbd}  Using the previous notations,
\[
\| Z_1 T_\lambda f\|_2 \less \lambda^{\frac12}
\|f\|_{L^{2,-\sigma}}, \qquad \| S_\lambda Z_1^* f\|_{L^{2,-\sigma}}
\less \lambda^{-\frac12} \|f\|_2
\]
for all $\lambda>1$.
\end{lemma}
\begin{proof}
First, \begin{equation}\label{eq:Z1T} Z_1 T_\lambda = w \la
\nabla\ra^{\frac12} T_\lambda + [\la \nabla\ra^{\frac12},w]
T_\lambda \end{equation}
 Now, by the same Fourier argument as above,
\[
\| \la \nabla\ra^{\frac12} T_\lambda f\|_{L^{2,-\sigma}} \less
\lambda^{\frac12} \|f\|_{L^{2,-\sigma}}
\]
Hence, the first term on the right-hand side of \eqref{eq:Z1T}
satisfies the desired bound. On the other hand, the commutator term
in~\eqref{eq:Z1T} can be written as
\[
\|[\la \nabla\ra^{\frac12},w] T_\lambda\|_{L^{2,-\sigma}\to L^2} \le
\| [\la \nabla\ra^{\frac12},w]w^{-1}\|_{L^2\to L^2} \| w
T_\lambda\|_{L^{2,-\sigma}\to L^2} \less 1
\]
uniformly in $\lambda$. Indeed, $[\la \nabla\ra^{\frac12},w]w^{-1}$
is a pseudo-differential operator of order zero and is therefore
$L^2$ bounded, whereas
\[
\| w T_\lambda\|_{L^{2,-\sigma}\to L^2} \less 1
\]
by the preceding. Next, we claim that
\begin{equation}\label{eq:claim4}\| Z_1 S_\lambda^* f\|_{2} \less
\lambda^{-\frac12} \|f\|_{L^{2,\sigma}}
\end{equation}
which will finish the proof by duality. To prove~\eqref{eq:claim4},
we write
\[
Z_1 S_\lambda^* = Z_1 T_\lambda T_\lambda^{-2} R_0(\lambda^2-i0)
\]
From \eqref{eq:algebra},
\[
\|T_\lambda^{-2} R_0(\lambda^2-i0) f\|_{L^{2,-\sigma}}\less
\lambda^{-1} \|f\|_{L^{2,\sigma}}
\]
provided $\sigma>\frac12$. Secondly, we have already shown that
\[Z_1 T_\lambda\::\:L^{2,-\sigma}\to L^2\] with bound
$\lambda^{\frac12}$. Thus, \eqref{eq:claim4} follows and we are
done.
\end{proof}

Now we continue with  the proof of \eqref{eq:largelam}. By the
resolvent identity, we have
\[
Z_1 R_L(\lambda^2+i0) Z_1^*=Z_1T_\lambda(I+S_\lambda L
T_\lambda)^{-1}S_\lambda Z_1^*
\]
provided $I+S_\lambda L T_\lambda$ is invertible as an operator on
$L^{2,-\sigma}$. This invertibility will follow by means of a
partial Neumann series via the following lemma. The proof of this
lemma, which is the crucial technical ingredient in this paper, will
be given in the next section.

\begin{lemma}
  \label{lem:largem} Given $A$ and $V$ as in Theorem~\ref{thm:main} as well as
  a positive constant $c>0$,
there exist  sufficiently large $m=m(c)$ and
$\lambda_0=\lambda_0(c)$ such that
\begin{equation}
  \label{eq:smallprod}
  \sup_{\lambda>\lambda_0}\| (R_0(\lambda^2+i0) L)^m \|_{L^{2,-\sigma}\to L^{2,-\sigma}}
  \le c
\end{equation}
Here $\sigma>4$.
\end{lemma}

In view of Lemmas~\ref{lem:STbd}, the estimate in
\eqref{eq:largelam} follows from the following result:

\begin{cor} \label{cor:main2} With the
notation from above  and for $\sigma>4$,  we have
$$
(I+S_\lambda L T_\lambda)^{-1}\::\:L^{2,-\sigma}\to L^{2,-\sigma}
$$
with a uniform norm for all large $\lambda$.
\end{cor}
\begin{proof}
We write the partial Neumann series, with $m$ as in
Lemma~\ref{lem:largem},
\[
(I+S_\lambda L T_\lambda)^{-1} = \Big(\sum_{k=0}^{m} (-1)^k
(S_\lambda L T_\lambda)^k \Big) (I+(-1)^{m+1}(S_\lambda L
T_\lambda)^{m+1})^{-1}
\]
By Lemma~\ref{lem:largem}, the inverse on the right-hand side exists
on $L^{2,-\sigma}$ with a uniform bound for all $\lambda>\lambda_0$.
Indeed, one has
\[
(S_\lambda L T_\lambda)^{m+1} = S_\lambda L (R_0(\lambda^2+i0)L)^m
T_\lambda
\]
so that, with some constant $C_1$ that only depends on $A$ and $V$,
\begin{align*}
&\|(S_\lambda L T_\lambda)^{m+1}\|_{L^{2,-\sigma}\to
L^{2,-\sigma}}\\
&\le  \|S_\lambda L\|_{L^{2,-\sigma}\to L^{2,-\sigma}} \|
(R_0(\lambda^2+i0)L)^m\|_{L^{2,-\sigma}\to L^{2,-\sigma}}\|
T_\lambda\|_{L^{2,-\sigma}\to L^{2,-\sigma}} \\
&\le C_1\, c <\frac12\end{align*} provided $c$ was chosen
sufficiently small.  Furthermore,
\[
S_\lambda L T_\lambda = 2i S_\lambda  A \cdot\nabla T_\lambda +
S_\lambda (i\div A + V) T_\lambda
\]
By \eqref{eq:Tlam} and Lemma~\ref{lem:Sagmon},
\[
\|S_\lambda  (i\div A + V) T_\lambda f\|_{L^{2,-\sigma}} \less
\|f\|_{L^{2,-\sigma}}
\]
Furthermore, again from \eqref{eq:Tlam} and Lemma~\ref{lem:Sagmon},
\[
\| S_\lambda  A \cdot\nabla T_\lambda \|_{L^{2,-\sigma}\to
L^{2,-\sigma}}\less \| S_\lambda  A\|_{L^{2,-\sigma}\to
L^{2,-\sigma}} \|\nabla T_\lambda \|_{L^{2,-\sigma}\to
L^{2,-\sigma}} \less \lambda^{-1} \lambda \less 1
\]
which means the finite sum of terms $k=0, \ldots, m$ can be controlled with
a bound independent of  $\lambda$.
\end{proof}

At this point the proof of Theorem~\ref{thm:limap} is essentially complete,
thanks to the identity
\begin{align*}
 \|\la\nabla\ra^\alpha &R_L(\lambda^2+i0)\la\nabla\ra^\alpha f\|_{L^{2,-\sigma}}
 =  \|\la\nabla\ra^\alpha T_\lambda (I + S_\lambda LT_\lambda)^{-1}
    S_\lambda \la\nabla\ra^{\alpha} f\|_{L^{2,-\sigma}} \\
 &\le \| \la\nabla\ra^\alpha T_\lambda\|_{L^{2,-\sigma}\to L^{2,-\sigma}}  
      \|(I + S_\lambda L T_\lambda)^{-1}\|_{L^{2,-\sigma}\to L^{2,-\sigma}}
      \|\la\nabla\ra^\alpha S_\lambda f\|_{L^{2,-\sigma}} \\
 &\les \la \lambda \ra^{2\alpha - 1}
 \|(I + S_\lambda L T_\lambda)^{-1}\|_{L^{2,-\sigma}\to L^{2,-\sigma}}
 \| f \|_{L^{2,\sigma}}
\end{align*}

For large $\lambda$, the desired operator bound for 
$(I + S_\lambda L T_\lambda)^{-1}$ is given by Corollary~\ref{cor:main2}.
For small $\lambda$, it follows from the Fredholm theory
arguments in Section~\ref{sec:small}.  One needs only to repeat the steps
taken in that section using the operator $T_\lambda^{-1}$ in place of
$Z_1$.

\section{The proof of Lemma~\ref{lem:largem}}

We start with the following observation: since $L=2i\nabla\cdot A -
i\div A + V$,
\begin{equation}
  \label{eq:teil}
  (R_0(\lambda^2+i0) L)^m = (2i)^m (R_0(\lambda^2+i0) \nabla \cdot A)^m
  + E_m(\lambda^2)
\end{equation}
where the error $E_m(\lambda^2)$ satisfies
\[
\| E_m(\lambda^2)\|_{L^{2,-\sigma}\to L^{2,-\sigma}} \le C(m,V,A)\,
\lambda^{-1}
\]
provided
\[
|A(x)|+|\div A(x)|+ |V(x)| \less \la x\ra^{-1-\eps}
\]
This follows from Agmon's limiting absorption
principle~\cite{Agmon}.

Thus, we are reduced to $L=\nabla\cdot A$. To deal with this case,
we shall perform a conical decomposition of the free resolvent. Let
$\{\chi_{\calS}\}_{\calS\in \Sigma}$ be a smooth partition of unity
on the sphere $S^2$ which is adapted to a family of caps $\Sigma$ of
diameter $\delta$ (which is a small parameter to be specified
later). For the most part, we shall drop the subscript $\calS$ so
that $\chi$ will denote any one of these cut-offs and $\chitil$ will
typically denote a cut-off associated to $\chi$ but with a dilated
cap as support. We write
\begin{equation}
  \label{eq:Rsplit}
  R_0(\lambda^2+i0)(x) = \sum_{\calS\in\Sigma} \frac{e^{i\lambda|x|}}{4\pi
|x|}\chi_\calS(x/|x|) =: \sum_{\calS\in\Sigma} R_\calS
(\lambda^2+i0)(x)
\end{equation}

We begin by studying the multiplier associated with $R_\calS$.

\begin{prop}
\label{prop:mult} Let $\chi$ be a cut-off supported in a
$\delta$-cap on $S^2$ where $\delta>0$ is a small parameter. Let
$K_\lambda$ be defined as
\[
K_\lambda(\xi):= \FT\Big[\frac{e^{i\lambda|x|}}{4\pi
|x|}\chi(x/|x|)\Big](\xi)
\]
where $\FT$ denotes the Fourier transform.  Then
\[
K_\lambda(\xi) := \left\{ \begin{array}{ll} O(\lambda^{-2}\delta^2)
& \wenn |\xi|<\frac{\lambda}{2} \\
O(|\xi|^{-2}) & \wenn |\xi|>10\lambda
\end{array} \right.
\]
and for $\frac{\lambda}{2}\le |\xi|\le 10\lambda$
\begin{equation}\label{eq:Klambda}
 K_\lambda(\xi) = O(\delta^{-2}\lambda^{-2}) +
\lambda^{-1}\chitil(\xi/|\xi|)f_\delta(\xi/\lambda) \big[
d\sigma_{\lambda S^2} (\xi) + i \PV \frac{1}{\lambda-|\xi|}\big]
\end{equation}
where $\chitil$ is a modified cut-off supported in twice the cap of
$\chi$ and $\|f_\delta\|_\infty \less 1$,
$\|f_\delta\|_{C^\alpha}\less \delta^{-2\alpha}$ for any $\alpha<1$.
\end{prop}
\begin{proof} By scaling, it suffices to set $\lambda=1$.  Let
\[
K(\xi) = K_{\eps,\delta}(\xi) = \int e^{-\eps|x|}
\frac{e^{i|x|}}{4\pi |x|} \chi(x/|x|) e^{-ix\cdot\xi}\, dx
\]
We assume that $\chi(x)$ is smooth and supported in a
$\delta$-neighborhood of $(0,0,1)$. Furthermore, by symmetry we can
assume that $\xi_2=0$. We shall use the identity
\begin{align}
K(\xi) &= \int_{S^2} \int_0^\infty e^{-\eps r} e^{ir} r \chi(\omega)
e^{-ir \omega\cdot\xi}\, dr d\sigma(\omega) \nn\\
&= \int_{S^2} (\eps -i(1-\omega\cdot\xi))^{-2} \chi(\omega)\,
d\sigma(\omega)\label{eq:K}
\end{align}

\smallskip
{\it Case 1:} $\xi_3\le \frac12$ and $|\xi|\le 10$.

\smallskip
Then, from \eqref{eq:K} we infer that
\[ K(\xi) = O(\delta^2) \]

\smallskip
{\it Case 2:} $|\xi_3|\ge \frac{|\xi|}{2}$ and $|\xi|> 10$.

\smallskip
In this case $|1-\omega\cdot\xi|\gtrsim |\xi|$ so that
\[ |K(\xi)|\less \frac{\delta^2}{|\xi|^2} \]
from \eqref{eq:K}.

\medskip

Cases 3 and 4 deal with $|\xi|>10$, $|\xi_3|<\frac{|\xi|}{2}$. Note
that then
\[
\{ \omega\cdot\xi \::\: \omega\in 2\calS \} = [a(\xi),b(\xi)]
\]
where $\calS:=\supp(\chi)\subset S^2$ and
$b(\xi)-a(\xi)\less\delta$. Moreover, $2\calS$ denotes the twice
dilated set $\calS$.

\smallskip
{\it Case 3:} $|\xi_3|\le \frac{|\xi|}{2}$ and $|\xi|> 10$, with $1\notin
[|\xi|a(\xi),|\xi|b(\xi)]$.

\smallskip Then
\begin{align*}
|K(\xi)| &\less \int_{a(\xi)+\delta}^{b(\xi)-\delta} \frac{\delta\,
ds}{(1-s|\xi|)^2} \less \frac{1}{|\xi|}
\int_{1-(b(\xi)-\delta)|\xi|}^{1-(a(\xi)+\delta)|\xi|}
\frac{\delta}{u^2}\, du \\
&\less \frac{\delta}{|\xi|} \big( |1-(b(\xi)-\delta)|\xi| |^{-1} +
|1-(a(\xi)+\delta)|\xi| |^{-1}\big) \\
&\less \frac{\delta}{|\xi|} \frac{1}{\delta|\xi|} \less |\xi|^{-2}
\end{align*}
as claimed

\smallskip
{\it Case 4:} $|\xi_3|\le \frac{|\xi|}{2}$ and $|\xi|> 10$, with $1\in
[|\xi|a(\xi),|\xi|b(\xi)]$.

\smallskip Here we write
\[
K(\xi) = \int_{I} \frac{\delta\psi(s)}{(s|\xi|-1-i\eps)^2}\, ds
\]
where $I$ is an interval of size $\sim\delta$ centered at
$|\xi|^{-1}$ and $|\psi^{(\ell)}(s)|\less \delta^{-\ell}$. Shifting
the center of $\psi$ to $0$ and abusing notation, we obtain
\begin{align*}
K(\xi) &= \int_{-c\delta}^{c\delta}
\frac{\delta\psi(s)}{(s|\xi|-i\eps)^2} \, ds =
\frac{\delta}{|\xi|}\int_{-c\delta}^{c\delta}\frac{\psi'(s)\, ds}{s|\xi|-i\eps}\\
&=
\frac{\delta}{|\xi|}\int_{-c\delta}^{c\delta}\frac{\psi'(s)-\psi'(0)}{s|\xi|-i\eps}
+ \frac{\delta}{|\xi|}\int_{-c\delta}^{c\delta}\frac{\psi'(0)\,
ds}{s|\xi|-i\eps} \\
& = O(|\xi|^{-2})
\end{align*}
using the bounds on $\psi'$ and $\psi''$.

\smallskip
{\it Case 5:} $\xi_3\ge\frac12$ and $\frac12 \le |\xi|\le 10$.

\smallskip In this case we write
\[
K(\xi) = O(\delta^{-2}) + \int_{\delta^{-2}}^\infty e^{-\eps r}
e^{ir} r a(r\xi)\, dr
\]
where
\[
a(r\xi) = \int_{S^2} \chi(\omega) e^{-ir\omega\cdot\xi}\,
d\sigma(\omega)
\]
By stationary phase
\[ a(r\xi) =\frac{e^{-ir|\xi|}}{r|\xi|} \Big(\chi(\xi/|\xi|) +
\chitil(\xi/|\xi|) \frac{\delta^{-2}}{|\xi|r}\Big) +
O\Big(\frac{\delta^{-4}}{|\xi|^3 r^3}\Big)
\]
Therefore, with $e:=\frac{\xi}{|\xi|}$,
\begin{align*}
K(\xi) &= O(\delta^{-2}) + \frac{\chi(e)}{|\xi|}
\frac{e^{[-\eps+i(1-|\xi|)]\delta^{-2}}}{\eps+i(1-|\xi|)} +
\frac{\chitil(e)}{|\xi|^2\delta^2} \int_{\delta^{-2}}^\infty
\frac{e^{[-\eps+i(1-|\xi|)]r}}{r}\, dr \\
&= O(\delta^{-2}) + \frac{1}{\eps-i(1-|\xi|)} \Big[
\frac{\chi(e)}{|\xi|} e^{[-\eps+i(1-|\xi|)]\delta^{-2}} +
\frac{\chitil(e)}{|\xi|^2} e^{[-\eps+i(1-|\xi|)]\delta^{-2}}\\
&\qquad\qquad  -\frac{\chitil(e)}{|\xi|^2\delta^2}
\int_{\delta^{-2}}^\infty
\frac{e^{[-\eps+i(1-|\xi|)]r}}{r^2}\, dr \Big] \\
&=: O(\delta^{-2}) + \frac{\chitil(e)}{\eps-i(1-|\xi|)}
f_{\eps,\delta}(\xi)
\end{align*}
Note that, as $\eps\to0$, $f_\delta:=\lim_{\eps\to0}
f_{\eps,\delta}$ satisfies
\[ \|f_\delta\|_\infty \less 1, \; \|f_\delta\|_{C^\alpha}\less
\delta^{-2\alpha} \] for any $\alpha<1$. Furthermore, in the sense
of distributions,
\[
\lim_{\eps\to0}\frac{\chitil(e)}{\eps-i(1-|\xi|)} = \chitil(e) \big[
d\sigma_{S^2} (\xi) + i \PV \frac{1}{1-|\xi|}\big]
\]
Here $\chitil$ on the right-hand side is modified to absorb any
needed constants.
\end{proof}

We shall use this result to prove Proposition~\ref{prop:limap}
below, which is a version of the limiting absorption principle.
First, we prove a lemma about the action of the singular part
in~\eqref{eq:Klambda} on functions.

\begin{lemma}
\label{lem:huh} Given a function $\phi$ in $\R^3$ and $0<\alpha<1$,
define
\[
[\phi]_\alpha(\xi) := \sup_{|h|<1}
\frac{|\phi(\xi)-\phi(\xi+h)|}{|h|^\alpha}
\]
Then
\[\begin{aligned}
&\Big|\int_{\R^3}\phi(\xi)\Big[ \sigma_{\lambda S^2}(d\xi)  +i\, \PV
\frac{d\xi}{\lambda-|\xi|} \chi_{[\lambda-1 < |\xi| < \lambda+1]}
\Big] \Big| \\
&\less \|\phi\|_{L^1(\lambda S^2)} + C_\alpha\,
\|[\phi]_\alpha\|_{L^1(\lambda S^2)}
\end{aligned}
\]
provided the right-hand side is finite.
\end{lemma}
\begin{proof}
It suffices to consider the principal value part. Thus,
\begin{align}
  &\Big| \PV \int_{||\xi|-\lambda|<1} \frac{\phi(\xi)}{|\xi|-\lambda}\,
  d\xi \Big| = \Big| \PV \int_{\lambda-1}^{\lambda+1} \frac{\beta^2
  \int_{S^2} \phi(\beta\theta) d\sigma(\theta) }{\beta-\lambda}\,d\beta
  \Big|\nn
  \\
  &\less \int_{\lambda-1}^{\lambda+1} \frac{\beta^2
  \int_{S^2} |\phi(\beta\theta)-\phi(\lambda\theta)| d\sigma(\theta)
  }{|\beta-\lambda|}\,d\beta\nn\\
&+ \Big| \PV \int_{\lambda-1}^{\lambda+1} \frac{\beta^2
  \int_{S^2} \phi(\lambda\theta) d\sigma(\theta) }{\beta-\lambda}\,d\beta
  \Big|\label{eq:2part}
\end{align}
The second term in \eqref{eq:2part} satisfies
\[\less \lambda \int_{S^2} |\phi(\lambda\theta)|\,d\sigma(\theta) \less
\lambda^{-1} \|\phi\|_{L^1(\lambda S^2)} \] whereas the first term
is
\[
\less \int_{\lambda-1}^{\lambda+1} \beta^2
|\beta-\lambda|^{\alpha-1} [\phi]_\alpha(\lambda\theta)\,
d\sigma(\theta) d\beta \le C_\alpha\, \|[\phi]_\alpha\|_{L^1(\lambda
S^2)}
\]
as claimed.
\end{proof}

We now turn to the limiting absorption principle. Note the decay
$\lambda^{-1}$ on the right-hand side which corresponds to a gain of
a derivative on the left-hand side. Also, note that the constant
does not depend on $\delta$ at least if $\lambda>\delta^{-2} $.

\begin{prop}\label{prop:limap}
Let  $ w=\la x\ra^{-\sigma}$ with $\sigma>4$. For
$\lambda>\delta^{-2}$ define the kernels
\begin{align*}
\tilde Q_\lambda(x,y) &:= w(x)
\frac{e^{i\lambda|x-y|}}{|x-y|}\chi\Big(\frac{x-y}{|x-y|} \Big) w(y)
\\
 Q_\lambda(x,y) &:= w(x)
\nabla_x\frac{e^{i\lambda|x-y|}}{|x-y|}\chi\Big(\frac{x-y}{|x-y|}
\Big) w(y)
\end{align*}
 Then,
\begin{align*}
&\|\tilde Q_\lambda\|_{2\to2} \le C_0 \lambda^{-1},\qquad
\|Q_\lambda\|_{2\to2} \le C_0
\end{align*}
The constant $C_0$ does not depend on~$\delta$.
\end{prop}
\begin{proof}
It will suffice to treat $Q_\lambda$.  We apply Schur's lemma. Thus,
using the notation of Proposition~\ref{prop:mult} (and assuming that
$w$ is real-valued)
\begin{align*}
&\int Q_\lambda(x,y) \bar{f}(y) g(x)\, dx dy \\&= \int \xi
K_\lambda(\xi)
\overline{\hat{w}\ast\hat{f}}(\xi) \hat{w}\ast\hat{g}(\xi)\, d\xi \\
&= \int\int \xi K_\lambda(\xi)
\overline{\hat{w}}(\xi-\xi_1)\hat{w}(\xi-\xi_2)\,d\xi
\overline{\hat{f}}(\xi_1)\hat{g}(\xi_2)\, d\xi_1d\xi_2
\end{align*}
The theorem follows provided we can show that
\begin{equation}
\label{eq:kernel} \sup_{\xi_2} \int\Big|  \int \xi K_\lambda(\xi)
\hat{w}(\xi_1-\xi)\hat{w}(\xi-\xi_2)\, d\xi \Big|\, d\xi_1
\less 1 %\lambda^{-1}
\end{equation}
First, note the bounds \begin{equation}\label{eq:hatw}
|\hat{w}(\xi)|\less \la \xi\ra^{-3-\eps},\qquad  |\nabla
\hat{w}(\xi)|\less \la \xi\ra^{-3-\eps}\end{equation} In fact, one
has rapid decay here but it is not needed. Second, it follows from
Proposition~\ref{prop:mult} that $K_\lambda:=
 K_1+K_2+K_3$ where \begin{equation} \label{eq:splitK}
\begin{aligned}
K_1(\xi) &= O(\delta^{-2}\lambda^{-2}) \chi_{[|\xi|<10\lambda]} \\
K_2(\xi) &= O(|\xi|^{-2})\chi_{[|\xi|>10\lambda]} \\
K_3(\xi) &= \lambda^{-1} \chi(e) f_\delta(\xi/\lambda)\Big[
d\sigma_{\lambda S^2}(\xi)  +i \PV \frac{1}{\lambda-|\xi|}
\chi_{[\lambda-1 < |\xi| < \lambda+1]} \Big]
\end{aligned}
\end{equation}
The cut-offs here are understood to be smooth. It is easy to see
that $K_1$ and $K_2$ contribute $O(\delta^{-2}\lambda^{-1})$ and
$O(\lambda^{-1})$ to~\eqref{eq:kernel}, respectively. To bound the
contribution of $K_3$, we use Lemma~\ref{lem:huh}. Thus, define
\[
\phi(\xi) := \xi \chi(\xi/|\xi|) f_\delta(\xi/\lambda)
\hat{w}(\xi_1-\xi)\hat{w}(\xi-\xi_2)
\]
Then
\begin{equation}\label{eq:J}
\|\phi\|_{L^1(\lambda S^2)} \less \lambda \int_{\lambda S^2}
\chi(\xi/|\xi|) \la\xi_1-\xi\ra^{-3-\eps}\la\xi-\xi_2\ra^{-3-\eps}
\, d\sigma(\xi) =: J_\lambda(\xi_1,\xi_2)
\end{equation}
as well as
\begin{equation}\label{eq:Jalpha}
\|[\phi]_\alpha \|_{L^1(\lambda S^2)} \less
\Big((\lambda\delta)^{-1}+(\delta^2\lambda)^{-\alpha}\Big)
J_\lambda(\xi_1,\xi_2) \less J_\lambda(\xi_1,\xi_2)
\end{equation}
provided $\lambda>\delta^{-2}$. In view of  Lemma~\ref{lem:huh} the
contribution by $K_3$ to~\eqref{eq:kernel} is bounded by
\[
\sup_{\xi_2} \lambda^{-1} \int J_\lambda(\xi_1,\xi_2) d\xi_1 \less
1%\lambda^{-1}
\]
and the proposition follows.
\end{proof}

Next, we study the effect of composing two resolvents which have
been restricted to disjoint conical regions.

\begin{prop}\label{prop:RL}
  Assume that $\sigma>4$ and
\begin{equation}
  \label{eq:decay}
  \sum_{|\alpha|\le 2} |D^\alpha
  \hat{A}(\xi)| \less \la \xi\ra^{-3-\eps} \qquad \forall \xi\in\R^3
\end{equation}
where $\eps>0$. Let $\calS_1, \calS_2\subset S^2$ with
$\dist(\calS_1,\calS_2)>5\delta$ where $\dist$ is the distance on
$S^2$. Let $R_1(\lambda^2)$ and $R_2(\lambda^2)$ be the free
resolvents which have been restricted to conical regions
corresponding to $\calS_1,\calS_2$, respectively. Then
\[
\| wR_{1}(\lambda^2) \nabla\cdot A R_{2}(\lambda^2) \nabla
w\|_{2\to2}\less \delta^{-2}\lambda^{-1}
\]
provided $\lambda>\delta^{-2}$.
\end{prop}
\begin{proof}
  We use Schur's lemma as in the proof of
  Proposition~\ref{prop:limap}. Thus, we write
\begin{align*}
 & \int\int\int g(x)w(x) \nabla_z R_1(\lambda^2)(x-z) A(z)\cdot \nabla_y R_2(\lambda^2)(z-y)
  w(y) \bar{f}(y)\, dxdydz \\
  & = \int\int \hat{g}(\xi) U(\xi,\eta)
  \bar{\hat{f}}(\eta)\,d\xi d\eta
\end{align*}
where (with real-valued $w$)
\[ U(\xi,\eta) := \int \hat{w}(\xi-\xi_1)
\xi_1 \widehat{R_1(\lambda^2)}(\xi_1) \hat{A}(\xi_2-\xi_1) \xi_2
\widehat{R_2(\lambda^2)}(\xi_2) \hat{w}(\eta-\xi_2)\, d\xi_1 d\xi_2
\]
We claim that
\begin{equation}\label{eq:schur}
  \sup_{\eta} \int_{\R^3} |U(\xi,\eta)|\, d\xi \less \delta^{-2}\lambda^{-1}
\end{equation}
By symmetry, this will imply the proposition. Next, we write as
in~\eqref{eq:splitK} for the Fourier transforms $K_\lambda^\pj=
\widehat{R_j(\lambda^2)}$ with $j=1,2$
\[
K_\lambda^\pj =  K_1^\pj+ K_2^\pj+ K_3^\pj
\]
The integral on the left-hand side of \eqref{eq:schur} is bounded by
\begin{equation} \sum_{i,j=1}^3 \int \Big|\int \hat{w}(\xi-\xi_1)
\xi_1 K_i^{\pone}(\xi_1) \hat{A}(\xi_2-\xi_1) \xi_2
K_j^{\ptwo}(\xi_2) \hat{w}(\eta-\xi_2)\, d\xi_1 d\xi_2 \Big|\, d\xi
\label{eq:nine}\end{equation} Of the nine different combinations
here all but $i=j=3$ are easy. Indeed,  if $i=1,2$ and for any
$j=1,2,3$,
\begin{align*}
&\int \Big|\int \hat{w}(\xi-\xi_1) \xi_1 K_i^{\pone}(\xi_1)
\hat{A}(\xi_2-\xi_1) \xi_2 K_j^{\ptwo}(\xi_2) \hat{w}(\eta-\xi_2)\,
d\xi_1
d\xi_2 \Big|\, d\xi \\
&\less \delta^{-2} \lambda^{-1} \int |\hat{w}(\eta-\xi_1)|\, d\xi_1
\int \Big|\int   \hat{A}(\xi_2-\xi_1) \xi_2 K_j^{\ptwo}(\xi_2)
\hat{w}(\eta-\xi_2)\,d\xi_2 \Big|\, d\xi d\xi_1 \\
& \less \delta^{-2}\lambda^{-1}
\end{align*}
by the discussion following \eqref{eq:kernel} (in particular,
recall~\eqref{eq:hatw}). It remains to consider $i=j=3$. For this we
shall use Lemma~\ref{lem:huh}. Let
\begin{align*}
G_\lambda(\xi_1,\eta) &:= \int \hat{A}(\xi_2-\xi_1)\xi_2
K_3^{\ptwo}(\xi_2) \hat{w}(\eta-\xi_2)\,  d\xi_2 \\
&= \lambda^{-1} \int \phi(\xi_2)\Big[ \sigma_{\lambda S^2}(d\xi_2)
+i\, \PV \frac{d\xi_2}{\lambda-|\xi_2|} \chi_{[\lambda-1 < |\xi_2| <
\lambda+1]} \Big]
\end{align*}
with
\[
\phi(\xi_2):=\hat{A}(\xi_2-\xi_1)\xi_2
\chi_2(\xi_2/|\xi_2|)f_\delta(\xi_2/\lambda) \hat{w}(\eta-\xi_2)
\]
Here $\chi_2$ is a cut-off adapted to $\calS_2$.  By
Lemma~\ref{lem:huh}, and \eqref{eq:J}, \eqref{eq:Jalpha},
\[
|G_\lambda(\xi_1,\eta)| \less  \int_{\lambda S^2}
\chi_2(\xi_2/|\xi_2|)
\la\xi_2-\xi_1\ra^{-3-\eps}\la\eta-\xi_2\ra^{-3-\eps} \,
d\sigma(\xi_2)
\]
Note that the same estimates hold if we replace $\hat{A} $ with
$\nabla\hat{A}$. Therefore,
\[
|\nabla_{\xi_1} G_\lambda(\xi_1,\eta)| \less \int_{\lambda S^2}
\chi_2(\xi_2/|\xi_2|)
\la\xi_2-\xi_1\ra^{-3-\eps}\la\eta-\xi_2\ra^{-3-\eps} \,
d\sigma(\xi_2)
\]
In view of these estimates we can apply Lemma~\ref{lem:huh} again to
obtain
\[
\begin{aligned}
&\Big|\int \hat{w}(\xi-\xi_1) \xi_1 K_3^{\pone} (\xi_1)
G_\lambda(\xi_1,\eta)\, d\xi_1 \Big|\\
&\less  \int_{\lambda S^2} \la \xi-\xi_1\ra^{-3-\eps}
\chi_1(\xi_1/|\xi_1|) \int_{\lambda S^2} \chi_2(\xi_2/|\xi_2|)
\la\xi_2-\xi_1\ra^{-3-\eps}\la\eta-\xi_2\ra^{-3-\eps} \,
d\sigma(\xi_2) \, d\sigma(\xi_1) \end{aligned}
\]
Hence the contribution of $i=j=3$ to \eqref{eq:nine} is bounded by
\begin{align*}
&  \int \int_{\lambda S^2} \int_{\lambda S^2}  \la
\xi-\xi_1\ra^{-3-\eps} \chi_1(\xi_1/|\xi_1|)\chi_2(\xi_2/|\xi_2|)
\la\xi_2-\xi_1\ra^{-3-\eps}\la\eta-\xi_2\ra^{-3-\eps} \,
d\sigma(\xi_2) d\sigma(\xi_1) \, d\xi \\
&\less   \int_{\lambda S^2} \int_{\lambda S^2}
\chi_1(\xi_1/|\xi_1|)\chi_2(\xi_2/|\xi_2|)
\la\xi_2-\xi_1\ra^{-3-\eps}\la\eta-\xi_2\ra^{-3-\eps} \,
d\sigma(\xi_2) d\sigma(\xi_1) \\
&\less \frac{1}{\lambda\dist(\calS_1,\calS_2)}\less
\lambda^{-1}\delta^{-1}.
\end{align*}
This is again smaller than $\delta^{-2} \lambda^{-1}$, as claimed.
\end{proof}

We now write the power on the right-hand side of~\eqref{eq:teil} as
a sum of products (dropping $\lambda^2+i0$ from the resolvent):

\begin{equation}
  \label{eq:sumprod}
  (R_0 \nabla \cdot A)^m =
  \sum_{\calS_1,\ldots,\calS_m\in\Sigma} R_{\calS_1} \nabla \cdot
  A \ldots \nabla \cdot
  A R_{\calS_m} \nabla \cdot A
\end{equation}
There are two types of chains $\calS_1,\calS_2, \ldots,\calS_m$ in
this sum:
\begin{itemize}
\item if $\dist(\calS_i,\calS_{i+1})\le 5\delta$ for all $1\le i\le
m-1$, then we call this chain {\em directed}
\item otherwise, we call it {\em undirected}
\end{itemize}

For the undirected chains there is the following corollary of the
previous proposition.

\begin{cor}
  \label{cor:undirect}
  If $\{\calS_j\}_{j=1}^m$ is undirected, then for $\sigma>4$
  \begin{equation}\label{eq:und}
\| R_{\calS_1} \nabla \cdot
  A \ldots \nabla \cdot
  A R_{\calS_m} \nabla \cdot A\|_{L^{2,-\sigma}\to
L^{2,-\sigma}} \le C(m,A) \delta^{-2}\lambda^{-1}
  \end{equation}
  provided $\lambda>\delta^{-2}$. In particular,
\begin{equation}\label{eq:und2}
\Big\|\sum_{\substack{\calS_1,\ldots,\calS_m\in\Sigma\\{\text{undirected}}}}
R_{\calS_1} \nabla \cdot
  A \ldots \nabla \cdot
  A R_{\calS_m} \nabla \cdot A \Big\|_{L^{2,-\sigma}\to L^{2,-\sigma}} \le C(m,A)
\delta^{-2(m+1)}\lambda^{-1}
\end{equation}
provided $\lambda>\delta^{-2}$.
\end{cor}
\begin{proof} This follows by applying Proposition~\ref{prop:RL} to
one pair of resolvents where $\dist(\calS_i,\calS_{i+1})> 5\delta$;
for the others, use Proposition~\ref{prop:limap}. More precisely,
with $i$ as specified, we write
\begin{equation}\label{eq:Aw}
AR_{\calS_i}\nabla\cdot A R_{\calS_{i+1}}\nabla\cdot A =
Aw^{-1}wR_{\calS_i}\nabla\cdot A R_{\calS_{i+1}}\nabla\cdot w w^{-1}
A
\end{equation}
where as usual $w(x)=\la x\ra^{-\sigma}$. In view of $|A(x)|\less
\la x\ra^{-2\sigma}$ and by our assumptions on $\hat{A}$, we apply
Proposition~\ref{prop:RL} to the right-hand side of~\eqref{eq:Aw} to
conclude that
\begin{equation}\label{eq:i}
\| wR_{\calS_i}\nabla\cdot A R_{\calS_{i+1}}\nabla\cdot w \|_{2\to2}
\less \delta^{-2}\lambda^{-1}
\end{equation}
To combine this with Proposition~\ref{prop:limap}, we insert factors
of $ww^{-1}$ as follows: with $\tilde A:= w^{-1} A w^{-1}$,
\begin{align*}
&  \prod_{j=1}^m (R_{\calS_j} \nabla A) =  w^{-1} \,(w R_{\calS_1}
\nabla w)\, \tilde A  \,(w R_{\calS_2} \nabla w)\, \tilde A \cdot\ldots\\
& \ldots \cdot \tilde A \,(wR_{\calS_i}\nabla\cdot A
R_{\calS_{i+1}}\nabla\cdot w)\, \tilde A  \,(w R_{\calS_{i+2}}
\nabla w)\, \ldots \,(w R_{\calS_{m}} \nabla w)\,\tilde A w
\end{align*}
Observe that
\[
\sup_j \| w R_{\calS_j} \nabla w \|_{2\to2} \le C
\]
uniformly in $\lambda>\delta^{-2}$ as well as $\|\tilde A f\|_2\less
\|f\|_2$. Combining this with \eqref{eq:i} yields~\eqref{eq:und}. To
pass to~\eqref{eq:und2} one sums over all possible choices of
undirected chains of which there are no more than $(C/\delta)^{2m}$.
\end{proof}
\begin{remark}
The summation over all possible paths is quite inefficient, as it does not
take advantage of any orthogonality between different operators
$R_\shell$.  However large the constants may be, once $A$, $m$, and $\delta$ are
fixed, the bound in \eqref{eq:und2} 
still approaches zero in the limit $\lambda \to \infty$.
\end{remark}

Finally, we turn to the directed chains. For these it will be
important that $\delta m\ll 1$ to ensure that the composition of
resolvents restricted to any directed chain remains outgoing.
Moreover, we will need to distinguish the {\em near} and {\em far}
parts of the free resolvent kernels which are defined as follows:
\begin{align*}
  Q_\calS^0 (x,y) &:= w(x) [\nabla_y R_\calS(x-y)] \chi(|x-y|<\rho)
  w (y) \\
  Q_\calS^1 (x,y) &:= w(x) [\nabla_y R_\calS(x-y)] \chi(|x-y|>\rho)
  w (y)
\end{align*}
where $1= \chi(|x-y|<\rho) + \chi(|x-y|>\rho)$ is a smooth partition
of unity adapted to the indicated sets. The parameter $\rho$ here is
a small number depending on~$m$. For the near part, we have the
following refinement of Proposition~\ref{prop:limap}.

\begin{prop}
  \label{prop:limap2}
  Under the conditions of Proposition~\ref{prop:limap} one has
  \[
\|  Q_\calS^0 \|_{2\to2} \le C_2 \rho,\qquad \| Q_\calS^1 \|_{2\to2}
\le C_2
  \]
  provided $\lambda>\delta^{-2}\rho^{-1}$. Here $C_2$ does not depend on
  $\delta$.
\end{prop}
\begin{proof}
Because of Proposition~\ref{prop:limap} it will suffice to prove the
bound on $Q_\calS^0$. In this proof, we shall write
\[
\chi_\rho(x-y):=\chi(|x-y|<\rho) \] Observe that $\hat\chi_\rho$ is
rapidly decaying outside of a ball of size $\less \rho^{-1}$. Thus,
as in the proof of Proposition~\ref{prop:limap}, and with $\tilde
K_\lambda(\xi):=\xi K_\lambda(\xi)$,
\begin{align*}
&\int Q_\calS^0(x,y) \bar{f}(y) g(x)\, dx dy \\&= \int [\tilde
K_\lambda\ast \hat{\chi_\rho}](\xi)
\overline{\hat{w}\ast\hat{f}}(\xi) \hat{w}\ast\hat{g}(\xi)\, d\xi \\
&= \int\int [\tilde K_\lambda\ast \hat{\chi_\rho}](\xi)
\overline{\hat{w}}(\xi-\xi_1)\hat{w}(\xi-\xi_2)\,d\xi
\overline{\hat{f}}(\xi_1)\hat{g}(\xi_2)\, d\xi_1d\xi_2
\end{align*}
The theorem follows provided we can show that
\begin{equation}
\label{eq:claim20} \sup_{\xi_2} \int\Big|  \int [\tilde
K_\lambda\ast \hat{\chi_\rho}](\xi)
\hat{w}(\xi_1-\xi)\hat{w}(\xi-\xi_2)\, d\xi \Big|\, d\xi_1 \less
\rho
\end{equation}
It follows from Proposition~\ref{prop:mult} that
\[ \tilde
K_\lambda:=
 \tilde K_1+\tilde K_2+\tilde K_3
 \]
 where (with smooth cut-offs)
\begin{align}
&[\tilde K_1\ast\hat{\chi_\rho}](\xi) = O(\delta^{-2}\lambda^{-1}) \label{eq:tilk1}\\
&[\tilde K_2\ast\hat{\chi_\rho}](\xi) = O(\lambda^{-1}) \label{eq:tilk2}\\
&\tilde K_3\ast\hat{\chi_\rho} =\nn \\&= \lambda^{-1}
\hat{\chi_\rho}\ast \Big\{\chi_\calS\, f_\delta(\cdot/\lambda)\Big[
\lambda d\sigma_{\lambda S^2}(\eta)  +i \PV
\frac{\eta}{\lambda-|\eta|} \chi_{[\lambda-1 < |\eta| < \lambda+1]}
\Big]\Big\} \label{eq:tilk3}
\end{align}
We also used there that $\lambda\gg\rho^{-1}$. The contributions of
\eqref{eq:tilk1} and~\eqref{eq:tilk2} to~\eqref{eq:claim20} are
treated as in Proposition~\ref{prop:limap} and yield a bound of
$\delta^{-2}\lambda^{-1}< \rho$ as desired. For the contribution
of~\eqref{eq:tilk3} we note that
\[
|\tilde K_3\ast\hat{\chi_\rho}|(\xi)\less \rho
\]
Hence, the contribution of \eqref{eq:tilk3} to \eqref{eq:claim20} is
controlled by
\begin{align*}
& \less \rho \sup_{\xi_2} \int \int
|\hat{w}(\xi_1-\xi)\hat{w}(\xi-\xi_2)|\, d\xi d\xi_1 \less \rho
\end{align*}
as desired.
\end{proof}

Next, we write
\begin{equation}\label{eq:big}\begin{aligned}
&\sum_{\substack{\calS_1,\ldots,\calS_m\in\Sigma\\{\text{directed}}}}
R_{\calS_1} \nabla \cdot
  A \ldots \nabla \cdot
  A R_{\calS_m} \nabla \cdot A \\
  &=
  \sum_{\substack{\calS_1,\ldots,\calS_m\in\Sigma\\{\text{directed}}}}\;\;
  \sum_{\eps_1,\ldots,\eps_m=0,1} w^{-1}\,
Q_{\calS_1}^{\eps_1}\,
  \tilde A \ldots \tilde
  A \,Q_{\calS_m}^{\eps_m}\, \tilde  A\, w
  \end{aligned}
\end{equation}
Fix a directed chain and assume without loss of generality that it
is directed along the positive $x_1$-axis. Since $\delta m\ll 1$,
one has \[ Q_{\calS_j}^1(x,y)=0 \text{\ \ unless\ \ } x_1-y_1>
\frac{\rho}{2}
\]
for each $1\le j\le m$. Next, we decompose
\[ \tilde A = \sum_{n\in\Z} \tilde A_n, \qquad \tilde A_n(x):=\tilde A(x)\chi_{[n\rho/2<x_1<(n+1)\rho/2]} \]
We start by estimating the contribution of products consisting
entirely of far kernels.

\begin{lemma}\label{lem:far} Suppose that $|A(x)|\le C_A \la x\ra^{-2\sigma-1-\eps}$ with $\sigma>4$.
Then, using the previous notations,
\[
\Big\| Q_{\calS_1}^{1}\,
  \tilde A \ldots \tilde
  A \,Q_{\calS_m}^{1}\, \tilde  A \Big\|_{2\to2} \le
  \frac{C_3^m}{m!\,\rho^m}
\]
provided $\lambda>\delta^{-2}+\rho^{-1}$. The constant $C_3$ here
depends only on $A$.
\end{lemma}
\begin{proof}
By our assumptions,
\[
\|\tilde A_n f\|_{2}\le  C_A (1+|n|\rho/2)^{-1-\eps} \|f\|_2
\]
Moreover, since $\sup_{1\le j\le m}\|Q_{\calS_j}^{1} \|_{2\to2}\le
C_2$,
\begin{align*}
  & \|Q_{\calS_1}^{1}\,
  \tilde A \ldots \tilde
  A \,Q_{\calS_m}^{1}\, \tilde  A \Big\|_{2\to2} \\
  &\le \sum_{n_1>n_2>\ldots>n_m} \|Q_{\calS_1}^{1}\,
  \tilde A_{n_1} \ldots \tilde
  A_{n_{m-1}} \,Q_{\calS_m}^{1}\, \tilde  A_{n_m} \Big\|_{2\to2}\\
  &\le C_2^m \sum_{n_1>n_2>\ldots>n_m} \prod_{j=1}^m \|\tilde
  A_{n_j}\|_{2\to2}\\
  &\le C_A^m C_2^m \sum_{n_1>n_2>\ldots>n_m} \prod_{j=1}^m
(1+|n_j|\rho/2)^{-1-\eps} \\
&\le \frac{C_A^m C_2^m }{m!}
\sum_{n_1,n_2,\ldots,n_m\in\Z}\prod_{j=1}^m
(1+|n_j|\rho/2)^{-1-\eps} \\
 &= \frac{C_3^m}{\rho^m m!}
\end{align*}
as claimed.
\end{proof}

 Next, we turn to the general case.

\begin{lemma}\label{lem:gen} Under the conditions of
Lemma~\ref{lem:far},
\[
\sum_{\eps_1,\ldots,\eps_m=0,1}\,\|Q_{\calS_1}^{\eps_1}\,
  \tilde A \ldots \tilde
  A \,Q_{\calS_m}^{\eps_m}\, \tilde  A\,\|_{2\to2} \le C_5^m
  m^{-\frac{m}{16}}
\]
where $C_5$ only depends on $A$.
\end{lemma}
\begin{proof} Let $\mu=\sum_{j=2}^m\eps_j$. Then
\begin{align}
&\sum_{\eps_1,\ldots,\eps_m=0,1}\,\|Q_{\calS_1}^{\eps_1}\,
  \tilde A \ldots \tilde
  A \,Q_{\calS_m}^{\eps_m}\, \tilde  A\,\|_{2\to2} \nn \\
&\le \sum_{\eps_1,\ldots,\eps_m=0,1}\;
{\sum_{n_1}}^{(\eps_2)}\ldots{\sum_{n_{m-1}}}^{(\eps_m)}\sum_{n_m}
C_2^m \rho^{1-\eps_1}\rho^{m-1-\mu} \prod_{j=1}^m \|\tilde
  A_{n_j}\|_{2\to2}  \label{eq:uff}
\end{align}
Here, for fixed $n_{i+1}$,
\[{\sum_{n_i}}^{(\eps_{i+1})} =\left\{
\begin{array}{cc}
\sum_{n_i>n_{i+1}} & \text{\ \ if\ \ }\eps_{i+1}=1 \\
\sum_{n_{i+1}+3\ge n_i\ge n_{i+1}} & \text{\ \ if\ \ }\eps_{i+1}=0
\end{array}\right.
\]
Now
\begin{align}
  \eqref{eq:uff} &\le  2\sum_{\eps_2,\ldots,\eps_m=0,1}\; {\sum_{n_1}}^{(\eps_2)}
  \ldots{\sum_{n_{m-1}}}^{(\eps_m)}\sum_{n_m} (C_A\,C_2)^m \nn\\
&\qquad\qquad\qquad\cdot \rho^{m-1-\mu} \prod_{j=1}^m (1+|n_j|\rho/2)^{-1-\eps} \nn \\
&\le  (4C_A\,C_2)^m \sum_{\ell=1}^{m-1} \binom{m-1}{\ell}
\frac{\rho^\ell}{(m-\ell-1)!}
\Big(\frac{C}{\rho}\Big)^{m-\ell-1}\label{eq:uff2}
\end{align}
by counting and symmetry as in the proof of Lemma~\ref{lem:far}.
Simplifying further, we conclude that
\begin{equation}
\label{eq:C4} \eqref{eq:uff2} \le C_4^m \rho^{-(m-1)}
\sum_{\ell=1}^{m-1} \binom{m-1}{\ell}
\frac{\rho^{2\ell}}{(m-\ell-1)!}
\end{equation}
The contribution of the sum over $\ell \ge
\frac{m-1}{2}+\frac{m}{4}$ to the right-hand side of~\eqref{eq:C4}
is at most $(2C_4)^m \rho^{\frac{m}{2}}$. On the other hand, the sum
over $\ell < \frac{m-1}{2}+\frac{m}{4}$ is bounded by
\[
(2C_4)^m\; \frac{\rho^{-(m-1)}}{\lfloor m/4\rfloor!}
\]
Setting $\rho:= m^{-\frac18}$ the lemma follows.
\end{proof}

Using \eqref{eq:big}, Lemma~\ref{lem:gen} and the observation that
there are at most
  $\delta^{-2}C^m$ directed chains we
conclude that
\begin{equation}\label{eq:d2}
\Big\|\sum_{\substack{\calS_1,\ldots,\calS_m\in\Sigma\\{\text{directed}}}}
R_{\calS_1} \nabla \cdot
  A \ldots \nabla \cdot
  A R_{\calS_m} \nabla \cdot A \Big\|_{L^{2,-\sigma}\to L^{2,-\sigma}} \le
  \delta^{-2}C_6^m
  m^{-\frac{m}{16}}
\end{equation}
Recall that in Lemma~\ref{lem:largem} we are given an operator $L$
(quickly reduced to the case $L = \nabla \cdot A$) and a small parameter 
$c > 0$.  Based on the value of $C_6(A)$ from \eqref{eq:d2} we choose
$m$ and $\delta = (10 m)^{-1}$ large enough so that the right side of 
\eqref{eq:d2} is less than $\frac{c}2$.  The bound for directed chains
is independent of $\lambda$.

For the undirected chains, we apply Corollary~\ref{cor:undirect}
directly.  With the quantities $m$ and $\delta$ already fixed, it is easy to
find $\lambda_0$ so that the right side of \eqref{eq:und2} is less than
$\frac{c}2$ whenever $\lambda > \lambda_0$.  This finishes the proof of
Lemma~\ref{lem:largem}.

%Combining this with the undirected case and choosing
%$\delta=m^{-2}$, say, finishes the proof of Lemma~\ref{lem:largem}.

\section{Appendix: absence of imbedded eigenvalues} \label{sec:absence}

 We consider $H=-\Delta+
\frac{i}{2}(A\cdot\nabla+\nabla\cdot A) +V$.

\begin{theorem}
\label{thm:froese} Assume that $V$ is bounded and converges to
zero at infinity and
\[|\nabla V(x)|, |A(x)|, |DA(x)| \le C \la x\ra^{-1-}.\]
Also assume that $\div A=0$. Then $H$ does not have any positive
eigenvalues.
\end{theorem}

Let $F\ge0$ be a radial, nondecreasing function with $|\nabla
F|\less 1$. Write $\nabla F=xg$ and let $\psi_F=e^F\psi$ for any
function $\psi$.  Suppose $H\psi = E\psi$ with $e^F\psi\in L^2$ and
$E>0$. We let $K $ be the symmetric generator of dilations:
\[ K = \frac12(x\cdot \nabla+\nabla\cdot x) \]
Then
\begin{align}
H\psi_F &= E\psi_F + [-(\nabla\cdot \nabla F+\nabla F\cdot
\nabla)+|\nabla F|^2+iA\cdot(\nabla F)] \psi_F.
\label{eq:HFcommute}
\\
\la \psi_F, H\psi_F\ra &= \la \psi_F,(|\nabla F|^2+E)\psi_F\ra. \label{eq:Eid}\\
\la \psi_F,[H,K]\psi_F\ra &= -4\|\sqrt{g}K\psi_F\|_2^2 +
\la\psi_F,C\psi_F\ra - 2 \Im \la A\cdot (\nabla F)\psi_F, K\psi_F\ra,\label{eq:HK1} \\
C &= (x\cdot \nabla)^2 g - x\cdot \nabla(|\nabla F|^2).
\label{eq:Cdef} \\
\la \psi_F,[H,K]\psi_F\ra &= 2E\|\psi_F\|_2^2 + \la \psi_F,(-i\tilde
A\cdot \nabla -2V-x\cdot\nabla V) \psi_F\ra, \label{eq:HK2} \\
\tilde A_j &= A_j + x\cdot \nabla A_j. \nn
\end{align}
Here $C$ is a multiplication operator, i.e., the derivatives only
act on the functions in the definition of~$C$. These are relatively
straightforward commutator identities. For example,  to derive
\eqref{eq:HK1} we proceed as follows:
\begin{align*}
 \la \psi_F,[H,K]\psi_F\ra &=  \la H\psi_F,K\psi_F\ra + \la
 K\psi_F,H\psi_F\ra
 = 2\Re \la H\psi_F,K\psi_F\ra
 \end{align*}
From \eqref{eq:HFcommute}, \begin{align*}
 H\psi_F &= E\psi_F +
[-g(\nabla\cdot x+ x\cdot
\nabla)+|\nabla F|^2+i(\nabla F)\cdot A] \psi_F - x\cdot(\nabla g) \psi_F \\
& = E\psi_F -2g K\psi_F + (|\nabla F|^2+i(\nabla F)\cdot A) \psi_F
- x\cdot(\nabla g)\psi_F
\end{align*}
Hence, \begin{equation} \label{eq:inter1}\begin{aligned} & \la
\psi_F,[H,K]\psi_F\ra \\
&= -4 \Re \la gK\psi_F, K\psi_F\ra + 2\Re \la (|\nabla F|^2 -
x\cdot(\nabla g)\psi_F,K\psi_F \ra \\
&\quad + 2 \Re \la i(A\cdot \nabla F)\psi_F,K\psi_F \ra
\end{aligned}
\end{equation}
Let $w$ be any real-valued function. Then, formally,
\begin{align*}
\la K\psi_F,w\psi_F\ra & = - \la \psi_F,K(w\psi_F)\ra \\
&= -\la\psi_F , [K,w]\psi_F\ra - \la\psi_F, wK\psi_F\ra \\
&= -\la \psi_F, [K,w] \psi_F\ra - \la w\psi_F, K\psi_F\ra
\end{align*}
and therefore,
\[ 2\Re\la K\psi_F,w\psi_F\ra = -\la\psi_F,[K,w]\psi_F\ra \]
Setting
\[ w = |\nabla F|^2- x\cdot(\nabla g)\] we can further simplify~\eqref{eq:inter1}
to~\eqref{eq:HK1}.

Let $\psi\in L^2$ and $E>0$ satisfy $H\psi=E\psi$. Let $\alpha>0$ be
a small constant. We will prove that $e^{\alpha|x|}\psi\in L^2$. To
this end define, for all $R>1$,
\[ F_R(r) = \alpha\int_0^r \chi(\rho)(1-\chi(\rho/R))\, d\rho \]
where $\chi(r)=0$ if $|r|<1$ and $\chi(r)=1$ if $|r|>2$, and
$\chi\ge0$ and smooth.
 Assume that $\|e^{F_R}\psi\|_2 \to\infty$ as
$R\to\infty$. Define $\phi_R= \psi_{F_R}/\|\psi_{F_R}\|_2$. Then
\[ \lim_{R\to \infty}\int_{|x|\le M} |\phi_R(x)|^2 \, dx =0 \]
for all $M>0$. In particular, \begin{equation}\label{eq:omega} \la
\phi_R, \omega \phi_R\ra \to 0\end{equation} as $R\to\infty$ for any
bounded $\omega$ with $|\omega(x)|\to0$ as $|x|\to\infty$.
By~\eqref{eq:Eid},
\begin{align*}
\|\nabla\phi_R\|_2^2 &= \la \phi_R, H\phi_R\ra + \la \phi_R,
(-iA\cdot\nabla -V)\phi_R\ra \\
&= \la \phi_R, (|\nabla F_R|^2 + E)\phi_R\ra + \la
\phi_R, (-iA\cdot\nabla -V)\phi_R\ra \\
&\le \|\nabla F_R\|_\infty^2 + E + \|V\|_\infty +
\|A\|_\infty^2/2 + \|\nabla\phi_R\|_2^2/2
\end{align*}
Since $\sup_{R>1}  \|\nabla F_R\|_\infty  <\infty$, it follows
that \begin{equation}\label{eq:sup_grad}\sup_{R>1} \|\nabla
\phi_R\|_2 <\infty. \end{equation} We now claim that
\begin{equation}\label{eq:claim1} \liminf_{R\to\infty} \la \phi_R,[H,K]\phi_R\ra \ge  2E
\end{equation}
This will lead to a contradiction via the second
identity~\eqref{eq:HK1} provided $\alpha$ is small depending on
$E>0$. To verify the claim, we need to check that
\[ \lim_{R\to\infty}  \la\phi_R ,(-i\tilde
A\cdot \nabla -2V-x\cdot\nabla V)\phi_R\ra =0,\] see \eqref{eq:HK2}.
However, this property follows immediately from~\eqref{eq:omega}
and~\eqref{eq:sup_grad} because of the decay of $\tilde A$
and~$V,x\cdot\nabla V$. Next, we use~\eqref{eq:HK1} to conclude that
\begin{align}
\label{eq:limsup} \limsup_{R\to\infty} \la \phi_R,[H,K]\phi_R\ra
& \le \limsup_{R\to\infty} |\la \phi_R,[(\x\cdot \nabla)^2 g_R]  \phi_R\ra  \\
&+\limsup_{R\to\infty} |\la \phi_R,[x\cdot \nabla(|\nabla F_R|^2)]\phi_R\ra \nn\\
&+2\limsup_{R\to\infty} |\la A\cdot (\nabla F_R)
\phi_R,K\phi_R\ra|.\nn
\end{align}
Note that
\begin{align}
\label{eq:limsupAK}\left|\la A\cdot (\nabla F_R) \phi_R,K\phi_R
\ra\right| &= \left|\la A\cdot (\nabla F_R) \phi_R,(3/2+x\cdot
\nabla) \phi_R  \ra\right|  \\
&\lesssim |\la A\cdot (\nabla F_R) \phi_R, \phi_R  \ra|+
\|A\cdot\nabla F_R |x|\phi_R\|_2\|\nabla \phi_R\|_2 \nn\\
&\to 0,\,\,\,\,\text{ as }R\to \infty. \nn
\end{align}
In the last line we used (\ref{eq:omega}), (\ref{eq:sup_grad}) and
 the decay of $|A|$ at infinity. Now,
\[ (r\partial_r)^2 g_R = \frac{\alpha}{r}
\chi(r)(1-\chi(r/R))-\alpha\chi'(r)+\frac{\alpha}{R}\chi'(r/R)+\alpha
r\chi''(r) - \frac{\alpha r}{R^2}\chi''(r/R) \] which implies that
\begin{equation}
\label{eq:bd1} \sup_{R>1} |(r\partial_r)^2 g_R(r)|\less
\frac{\alpha}{\la r\ra}.
\end{equation}
Thus, using (\ref{eq:omega}), we have
\begin{align}
\label{eq:limsupg}\lim_{R\to\infty}  \la \phi_R, [(\x\cdot
\nabla)^2 g_R ] \phi_R\ra = 0.
\end{align}
Finally,
\[
(r\partial_r)(\nabla F_R)^2 = 2\alpha^2 r\chi'(r)\chi(r)  -2
\frac{r\alpha^2}{R}(1-\chi(r/R))\chi'(r/R)
\]
which yields
\begin{equation}
\label{eq:bd2} \sup_{R>1} |(r\partial_r)(\nabla F_R)^2(r)|\less
\alpha^2.
\end{equation}
Using (\ref{eq:limsupAK}), (\ref{eq:limsupg}) and (\ref{eq:bd2})
in (\ref{eq:limsup}), we obtain
\[ \limsup_{R\to\infty} \la \phi_R,[H,K]\phi_R\ra  \less
\alpha^2.
\]
For small $\alpha\le \alpha_0(E)$ we obtain a contradiction to
\eqref{eq:claim1}.

Next, we claim that $e^{\alpha|x|}\psi\in L^2$ for all $\alpha>0$.
This can be done inductively, by increasing $\alpha$ in steps of
$\eps$ for suitable $\eps=\eps(\alpha)$. More precisely, with any
$\alpha>0$, we define
\begin{equation} \label{eq:FRdef} F_R(r) =
\alpha r \chi(r) + \eps\int_0^r \chi(\rho)(1-\chi(\rho/R))\,
d\rho.
\end{equation}
Then, on the one hand, \eqref{eq:claim1} remains unchanged.
 On the other hand, \eqref{eq:limsupAK} and \eqref{eq:limsupg} remain unchanged, and hence we
 have
\begin{align}
 \limsup_{R\to\infty} \la \phi_R,[H,K]\phi_R\ra &\le  \limsup_{R\to\infty} |\la \phi_R,  [(x\cdot \nabla) |\nabla
F_R|^2] \phi_R\ra|. \label{eq:bd4}
\end{align}
To bound the latter, we observe from~\eqref{eq:FRdef} that
\begin{align*}
 & r\partial_r|\nabla
F_R|^2 = r\partial_r \Big[\alpha\chi(r)+\alpha r\chi'(r) + \eps
\chi(r)(1-\chi(r/R)) \Big]^2 \\
&= 2r \Big[\alpha\chi(r)+\alpha r\chi'(r) + \eps
\chi(r)(1-\chi(r/R)) \Big] \\
&\qquad \cdot\Big[2\alpha\chi'(r)+\alpha r\chi''(r) + \eps \chi'(r)-
\eps R^{-1}\chi'(r/R) \Big]
\end{align*}
Thus,
\[ \Big|r\partial_r|\nabla
F_R|^2\Big| \less \alpha^2 \chi_{[|r|\le 2]} + (\alpha+\eps )\eps \]
whence
\[  \limsup_{R\to\infty} \la \phi_R,[H,K]\phi_R\ra  \less (\alpha+\eps
)\eps \] see \eqref{eq:bd4}. It follows that as long as $\eps\less
\alpha^{-1}$, this contradicts~\eqref{eq:claim1}. Since $
\int_1^\infty \alpha^{-1}\,d\alpha=\infty$, we see that
$e^{\alpha|x|} \psi\in L^2$ for all $\alpha>0$.

The final step in the proof of the theorem is the following lemma.

\begin{lemma}
\label{lem:UC} Let $H$ be as in the theorem. Assume that $\psi$
satisfies $H\psi=E\psi$ with $E>0$ and $e^{\alpha|x|}\psi\in L^2$
for all $\alpha>0$. Then $\psi=0$.
\end{lemma}
\begin{proof}
Let $F_\alpha = \alpha\la x\ra$ and $\psi_\alpha =
e^{F_\alpha}\psi$. Then
\begin{align*}
\|\nabla\psi_\alpha\|_2^2&=\la \psi_\alpha , -\Delta \psi_\alpha \ra
\ge \la \psi_\alpha, H \psi_\alpha\ra - C\|\psi_\alpha\|_2^2 -
 \|\nabla \psi_\alpha\|_2^2
\end{align*}
and therefore, by \eqref{eq:Eid},
\begin{align}
\label{eq:lowbnd}\|\nabla\psi_\alpha\|_2^2 &\ge \frac12\la
\psi_\alpha, H
\psi_\alpha\ra - C\|\psi_\alpha\|_2^2 \\
& \ge\frac12\la \psi_\alpha, |\nabla F_\alpha|^2 \psi_\alpha\ra -
C\|\psi_\alpha\|_2^2 \nn \\
&= \frac12\la \psi_\alpha, \alpha^2 r^2\la r\ra^{-2}
\psi_\alpha\ra - C\|\psi_\alpha\|_2^2.\nn
\end{align}
Since $[H,K]= -2\Delta - x\cdot\nabla V - i x_k (\partial_k
A_j)\partial_j $, we conclude that \begin{equation} \label{eq:upbnd}
\begin{aligned}
\|\nabla\psi_\alpha\|_2^2 &\le \la \psi_\alpha, [H,K] \psi_\alpha
\ra + C\|\psi_\alpha\|_2^2 \\
&\le \la \psi_\alpha, [(\x\cdot \nabla)^2 g_\alpha - x\cdot
\nabla(|\nabla F_\alpha|^2)]  \psi_\alpha \ra \\
& + |\la A\cdot(\nabla F_\alpha)\psi_\alpha, K\psi_\alpha\ra|  +
C\|\psi_\alpha\|_2^2.
\end{aligned}
\end{equation}
Note that
\begin{align}
  |\la A\cdot (\nabla F_\alpha)
\psi_\alpha,K\psi_\alpha \ra | &=  |\la A\cdot (\nabla F_\alpha)
\psi_\alpha,(3/2+x\cdot
\nabla) \psi_\alpha  \ra | \nn \\
&\le  C \|  \psi_\alpha\|_2^2 + \frac12 \|\nabla
\psi_\alpha\|_2^2. \nn
\end{align}
Using this in (\ref{eq:upbnd}), we obtain
\begin{align}
\|\nabla\psi_\alpha\|_2^2 &\le 2 \la \psi_\alpha, [(\x\cdot
\nabla)^2 g_\alpha - x\cdot \nabla(|\nabla F_\alpha|^2)] \psi_\alpha
\ra
+ C\|\psi_\alpha\|_2^2 \nn\\
&\le 2 \la \psi_\alpha, \big\{\alpha[3r^4\la r\ra^{-5}-2r^2\la
r\ra^{-3}]-2\alpha^2r^2\la r\ra^{-4} \big\}\psi_\alpha \ra   +
C\|\psi_\alpha\|_2^2. \label{eq:upbnd1}
\end{align}
Combining (\ref{eq:lowbnd}) and (\ref{eq:upbnd1}), we get
\[ \la \psi,e^{2F_\alpha}\big[\alpha^2(\frac12 r^2\la r\ra^{-2}+4r^2\la
r\ra^{-4}) + \alpha(4r^2\la r\ra^{-3}-6r^4\la r\ra^{-5})
-C\big]\psi\ra \le 0\] This can be written as
\[ \int |\psi(x)|^2 w_\alpha(x)\, dx \le 0\]
where $\inf_x w_\alpha(x)\ge c_0\alpha^2  >0 $ for all $\alpha$
large (with $c_0$ independent of $\alpha$). This is a contradiction.
\end{proof}


\begin{thebibliography}{99}

\bibitem{Agmon} Agmon, S. {\em Spectral properties of Schr\"odinger
operators and scattering theory.} Ann.\ Scuola Norm.\ Sup.\ Pisa Cl.\ Sci. (4) 2 (1975), no.~2, 151--218.

\bibitem{CK} Christ, M., Kiselev, A. {\em Maximal functions associated with
filtrations,} J.\ Funct.\ Anal.\ 179 (2001), 409-425.

\bibitem{CS1} Constantin, P., Saut, J.-C.
{\em Local smoothing properties of dispersive equations.} J.\ Amer.\
Math.\ Soc.~1 (1988), no.~2, 413--439.

\bibitem{CS2} Constantin, P., Saut, J.-C.
{\em Local smoothing properties of Schr\"odinger equations.} Indiana
Univ.\ Math.\ J.~38 (1989), no.~3, 791--810.

\bibitem{CFKS} Cycon, H. L., Froese, R. G., Kirsch, W., Simon, B.
{\em Schr\"odinger operators with application to quantum mechanics
and global geometry.}
 Texts and Monographs in Physics. Springer Study Edition. Springer-Verlag, Berlin, 1987.

\bibitem{FHOO} Froese, R., Herbst, I., Hoffmann-Ostenhof, M., Hoffmann-Ostenhof, T. {\em
 On the absence of positive eigenvalues for one-body Schr\"odinger operators.}  J.\ Analyse Math.~41  (1982), 272--284.

\bibitem{GST} Georgiev, V., Stefanov, A., Tarulli, M. {\em Smoothing
-  Strichartz estimates for the Schr\"odinger equation with small
magnetic potential}, preprint 2005.

\bibitem{GS} Goldberg, M., Schlag, W. {\em Dispersive estimates for Schr\"odinger operators in dimensions one and three.}
 Comm.\ Math.\ Phys.~251  (2004),  no.~1, 157--178.

\bibitem{IS} Ionescu, A., Schlag, W. {\em Agmon--Kato--Kuroda theorems for a large class of
perturbations.},  Duke Math.\ J.~131  (2006),  no.~3, 397--440.

\bibitem{JenKat} Jensen, A., Kato, T. {\em Spectral properties of Schr\"odinger operators and time-decay of
the wave functions.} Duke Math.\ J.~46  (1979), no.~3, 583--611.

\bibitem{JSS}  Journ\'e, J.-L.,  Soffer, A., Sogge, C.\ D. {\em Decay estimates for Schr\"odinger operators.}
Comm.\ Pure Appl.\ Math.~44 (1991), no.~5, 573--604.

\bibitem{Kato} Kato, T. {\em Wave operators and similarity for some non-selfadjoint operators.} Math.\ Ann.~162
(1965/1966),  258--279.

\bibitem{KatoY} Kato, T., Yajima, K. {\em Some examples of smooth
operators and the associated smoothing effect.} Rev.\ Math.\ Phys.~1
(1989), no.~4, 481--496.

\bibitem{KT} Keel, M., Tao, T. {\em Endpoint Strichartz estimates.}  Amer.\ J.\ Math.~120  (1998),  no.~5, 955--980.

\bibitem{KochT} Koch, H., Tataru, D. {\em Carleman estimates and absence of embedded eigenvalues}, preprint 2005.

\bibitem{rauch} Rauch, J. {\em Local decay of scattering solutions to Schr\"odinger's equation.} Comm.\ Math.\ Phys.~61  (1978), no.~2, 149--168.

\bibitem{RS4} Reed, M., Simon, B. {\em Methods of modern mathematical physics. IV. Analysis of operators.}
 Academic Press [Harcourt Brace Jovanovich, Publishers], New York-London, 1978.

\bibitem{Rob} Robert, D. {\em Asymptotique de la phase de diffusion \`a haute
\'energie pour des perturbations du second ordre du Laplacien}. (French)
Ann.\ Sci.\ \'Ec.\ Norm.\ Sup., IV S\'er.\ 25 (1992), No. 2, 107--134.

\bibitem{RS} Rodnianski, I., Schlag, W. {\em Time decay for solutions of
Schr\"odinger equations with rough and time-dependent potentials.}
Invent.\ Math.\ 155 (2004), no. 3, 451--513.

\bibitem{Sch} Schlag, W., {\em Dispersive estimates for Schr\"odinger operators: A survey}, preprint 2005.

\bibitem{SJ} Sj\'olin, P. {\em Regularity of solutions to Schr\"odinger equations.} Duke Math.\ J.~55 (1987),
699--715.

\bibitem{Stef} Stefanov, A. {\em Strichartz estimates for the magnetic
Schr\"odinger equation}, preprint 2004.

\bibitem{strich} Strichartz, R. {\em Restrictions of Fourier transforms to quadratic
surfaces and decay of solutions of wave equations.}  Duke Math.\
J.~44  (1977), no.~3, 705--714.

\bibitem{VE} Vega, L. {\em Schr\"odinger equations: Pointwise
convergence to the initial data.} Proc.\ Amer.\ Math.\ Soc.~102
(1988), 874--878.

\end{thebibliography}
\end{document}